\documentclass[12pt, reqno]{amsart}
\usepackage{latexsym}
\usepackage{amsfonts}
\usepackage{amssymb}
\usepackage{color}
\usepackage{graphicx}
\usepackage[mathscr]{euscript}

\newtheorem{theorem}{Theorem}[section]
\newtheorem{lemma}[theorem]{Lemma}
\newtheorem{proposition}[theorem]{Proposition}
\newtheorem{corollary}[theorem]{Corollary}
\theoremstyle{definition}
\newtheorem{definition}{Definition}[section]
\newtheorem{example}{Example}[section]

\theoremstyle{remark}
\newtheorem{remark}{Remark}[section]
\numberwithin{equation}{section}

\newcommand{\CC}{\mathbb{C}}

\newcommand{\NN}{\mathbb{N}}

\newcommand{\MM}{\mathbb{M}}

\newcommand{\supp}{\mathrm{supp}}
\newcommand{\diag}{\mathrm{diag}}

\newcommand{\rank}{\mathrm{rank}}

\newcommand{\im}{\mathrm{Im}}
\newcommand{\tr}{\mathrm{tr}}

\newcommand{\TT}{\mathscr{T}}
\newcommand{\FF}{\mathscr{F}}

\newcommand{\Lgl}{\mathfrak{gl}}

\begin{document}
%%%%%%%%%%%%%%%%%%%%%%%%%%%%%%%%%%%%%%%%%%%
\title[Nonabelian harmonic analysis and Functional Equations]
{Nonabelian harmonic analysis and functional equations on compact
groups}
\thanks{}
\author[J. An {\protect \and} D. Yang]
{Jinpeng An {\protect
\and} Dilian Yang}
\address{JINPENG AN, School of Mathematical Sciences, Peking University, Beijing, 100871, China}
\email{anjinpeng@gmail.com}
\address{DILIAN YANG,
Department of Mathematics $\&$ Statistics, University of Windsor, Windsor, ON N9B 3P4, CANADA}
\email{dyang@uwindsor.ca}

\begin{abstract}
Making use of nonabelian harmonic analysis and representation
theory, we solve the functional equation
$$f_1(xy)+f_2(yx)+f_3(xy^{-1})+f_4(y^{-1}x)=f_5(x)f_6(y)$$ on
arbitrary compact groups. The structure of its general solution is
completely described. Consequently, several special cases of the
above equation, in particular, the Wilson equation and the d'Alembert long equation,
are solved on compact groups.
\end{abstract}



\thanks{
Part of this work was conducted when the authors were at the
University of Waterloo.\\ \indent\hskip .1cm The second author was
partially supported by an NSERC grant. }

\date{}
\maketitle

%%%%%%%%%%%%%%%%%%%%%%%%%%%%%%%%%%%%%%%%%%%%%%
\section{Introduction}

Let $G$ be a group. The \textit{d'Alembert equation}
\begin{eqnarray}
\label{E:D'Alembert} f(xy)+f(xy^{-1})=2f(x)f(y),
\end{eqnarray}
where $f:G\to \CC$ is the function to determine, has a long history
(see \cite{Aczel}). It is easy to check that if $\varphi$ is a
homomorphism from $G$ into the multiplicative group of nonzero
complex numbers, the function $f(x)=(\varphi(x)+\varphi(x)^{-1})/2$
is a solution of Eq.~\eqref{E:D'Alembert} on $G$. Such solutions and
the zero solution are called classical solutions. Kannappan
\cite{Kan} proved that if $G$ is abelian, then all solution of
Eq.~\eqref{E:D'Alembert} are classical. This was generalized to
certain nilpotent groups in \cite{Cor1, Cor, Fri, Stet1, Stet2}. On
the other hand, Corovei \cite{Cor1} constructed a nonclassical
solution of Eq.~\eqref{E:D'Alembert} on the quaternion group $Q_8$.
It was realized later that Corovei's solution is nothing but the
restriction to $Q_8$ of the normalized trace function $\tr/2$ on
$SU(2)$, which is a nonclassical solution of
Eq.~\eqref{E:D'Alembert} on $SU(2)$ (c.f. \cite{ACN,Yang06}).
Recently, it was proved in \cite{Yang06,Yangthesis} that any
nonclassical continuous solution of Eq.~\eqref{E:D'Alembert} on a
connected compact group factors through $SU(2)$, and that the
function $\tr/2$ is the only nonclassical continuous solution on
$SU(2)$. This was generalized by Davison to arbitrary compact groups
in \cite{Davison1}, and further to any topological groups in
\cite{Davison2} (with the group $SU(2)$ replaced by $SL(2,\CC)$).
Hence Eq.~\eqref{E:D'Alembert} on topological groups has been
completely solved. For more results related to
Eq.~\eqref{E:D'Alembert}, we refer to \cite{Cho, Cho1, Davison3,
Stet4, Stet5, Stet6} and the survey \cite{Stet3}.

A well-known generalization of the d'Alembert equation is the
\textit{Wilson equation}
\begin{equation}\label{E:Wilson}
f(xy)+f(xy^{-1})=2f(x)g(y),
\end{equation}
where $f$ and $g$ are unknown complex functions on $G$. It was first
considered by Wilson \cite{Wil} and has also been extensively
studied (see \cite{Davison2, Davison3, Fri, Stet2, Stet3} and the
references therein). It turns out in \cite{Davison2} that
Eq.~\eqref{E:Wilson} is directly related to
Eq.~\eqref{E:D'Alembert}, where solutions of Eq.~\eqref{E:Wilson}
were used to construct the homomorphism $G\rightarrow SL(2,\CC)$
mentioned in the previous paragraph. Furthermore, it was shown (see,
e.g., \cite{Stet2}) that if $f$ and $g$ satisfy Eq.~\eqref{E:Wilson}
and $f\not\equiv0$, then $g$ is a solution of the \textit{d'Alembert
long equation}
\begin{equation}\label{E:long}
f(xy)+f(yx)+f(xy^{-1})+f(y^{-1}x)=4f(x)f(y).
\end{equation}
The question of solving Eq.~\eqref{E:long} on arbitrary topological
groups was raised in \cite{Davison1}. However, the approaches in
\cite{Davison1, Davison2, Yang06} do not apply to
Eqs.~\eqref{E:Wilson} and \eqref{E:long}.

The purpose of this paper is to study the equation
\begin{equation}\label{E:eq-general}
f_1(xy)+f_2(yx)+f_3(xy^{-1})+f_4(y^{-1}x)=f_5(x)f_6(y),
\end{equation}
where $f_i:G\to \CC$ ($i=1,...,6$) are unknown functions. It is
clear that Eq.~\eqref{E:eq-general} includes
Eqs.~\eqref{E:D'Alembert}--\eqref{E:long} as special cases. We will
find all $L^2$-solutions of Eq.~\eqref{E:eq-general} on arbitrary
compact groups. Consequently, we will solve Eqs.~\eqref{E:Wilson}
and \eqref{E:long} on compact groups completely. Here, it is worth
mentioning that, under some mild conditions, nonzero $L^2$-solutions
of the d'Alembert equation \eqref{E:D'Alembert} (or some of its
variant forms) exist only when $G$ is compact (c.f. \cite{PenRuk}).

Our main ingredients are nonabelian harmonic analysis on compact
groups and representation theory. Let $G$ be a compact group. Then
the Fourier transform transforms a square integrable function $f$ on
$G$ into an operator-valued function $\hat{f}$ on $\hat{G}$, the
unitary dual of $G$. Applying the Fourier transform to both sides of
Eq.~\eqref{E:eq-general} and taking some representation theory into
account, we will convert Eq.~\eqref{E:eq-general} into a family of
matrix equations. We call a tuple of matrices satisfying such matrix
equations an admissible (matrix) tuple. There are three types of
admissible tuples, i.e., complex, real, and quaternionic types,
which correspond to the three types of the representations
$[\pi]\in\hat{G}$, respectively. To determine the admissible tuples
is a question of linear algebra. We will find all admissible tuples
of each type. Then applying the Fourier inversion formula, we obtain
the general solution of Eq.~\eqref{E:eq-general}.

The structure of the general solution of Eq.~\eqref{E:eq-general}
can be compared with that of linear differential equations, where
any solution is the sum of a particular solution and a solution of
the associated homogeneous differential equation. In our case, the
\textit{homogeneous equation} associated with
Eq.~\eqref{E:eq-general} is
\begin{equation}\label{E:eq-homogeneous}
f_1(xy)+f_2(yx)+f_3(xy^{-1})+f_4(y^{-1}x)=0.
\end{equation}
It is obvious that the solutions of Eq.~\eqref{E:eq-homogeneous}
form a closed subspace of $L^2(G)^4$, and that the sum of a solution
of Eq.~\eqref{E:eq-general} and a solution of
Eq.~\eqref{E:eq-homogeneous} is also a solution of
Eq.~\eqref{E:eq-general}. Some obvious solutions of
Eq.~\eqref{E:eq-homogeneous} are provided by central functions. We
will determine the orthogonal complement of these obvious solutions
in the solution space of Eq.~\eqref{E:eq-homogeneous} by
constructing a spanning set using irreducible representations of $G$
into $O(1)$, $O(2)$, and $SU(2)$. We will also prove that any
solution of Eq.~\eqref{E:eq-general} is the sum of a solution of
Eq.~\eqref{E:eq-homogeneous} and a pure normalized solution of
Eq.~\eqref{E:eq-general} (see Section 3 for the definitions), and
will determine all pure normalized solutions of
Eq.~\eqref{E:eq-general}, which correspond to irreducible
representations of $G$ into $U(1)$, $O(2)$, $SU(2)$, and $O(3)$.
This provides a complete picture of the general solution of
Eq.~\eqref{E:eq-general}. These results will be proved in Theorems
\ref{T:nonhomogeneouspure}--\ref{T:nonhomogeneous}. As applications,
we will solve several special cases of Eq.~\eqref{E:eq-general},
including Eqs.~\eqref{E:Wilson} and \eqref{E:long}. In particular,
we will show that all nontrivial solutions of
Eqs.~\eqref{E:Wilson} and \eqref{E:long} factor through $SU(2)$, and
that the general solutions of Eq.~\eqref{E:long}
and Eq.~\eqref{E:D'Alembert} are the same.

The paper is organized as follows. Some basic properties of the
Fourier transform on compact groups and some facts in representation
theory will be briefly reviewed in Section 2. In Section 3 we will
give some basic definitions related to Eq.~\eqref{E:eq-general},
introduce the notion of admissible matrix tuples, reveal their
relations with Eq.~\eqref{E:eq-general}, and present some examples
which are the building blocks of the general solution. Then in
Section 4 we will determine all admissible matrix tuples. The main
results will be proved in Section 5. The general solutions of
several special cases of Eq.~\eqref{E:eq-general} will be given in
Section 6.

We should point out that one could apply our method in this paper to
some other types of functional equations on compact groups, and that
the method may be also generalized to solve functional equations on
non-compact groups admitting Fourier transforms.

Throughout this paper, $G$ denotes a compact group, $dx$ the
normalized Haar measure on $G$, and $L^2(G)$ the Hilbert space of
all square integrable functions on $G$ with respect to $dx$. By
solutions of Eq.~\eqref{E:eq-general} (or its special cases) on $G$
we always mean its $L^2$-solutions.

We would like to thank Professor H. Stetk{\ae}r for giving many
valuable comments.

%%%%%%%%%%%%%%%%%%%%%%%%%%%%%%%%%%%%%%%%%%%%%%

\section{Preliminaries}

As mentioned in the introduction, our basic tools in this paper are
Fourier analysis on compact groups and some results in
representation theory. In this section, we briefly review some
fundamental facts in these two subjects that will be used later.

\subsection{Fourier analysis}

We mainly follow the approach of \cite[Chapter 5]{Fo}. Let $\hat{G}$
be the unitary dual of the compact group $G$. For $[\pi]\in\hat{G}$,
we view $\pi$ as a homomorphism $\pi:G\rightarrow U(d_\pi)$, where
$d_\pi$ is the dimension of the representation space. Let
$\MM(n,\CC)$ denote the space of all $n\times n$ complex matrices.
For $f\in L^2(G)$, the Fourier transform of $f$ is defined by
$$\hat{f}(\pi)=d_\pi\int_Gf(x)\pi(x)^{-1}dx\in\MM(d_\pi,\CC), \quad [\pi]\in\hat{G}.$$
Note that for the sake of convenience, our definition is different
from the one in \cite{Fo} by a factor $d_\pi$. In our setting, the
Fourier inversion formula is
$$f(x)=\sum_{[\pi]\in\hat{G}}\tr(\hat{f}(\pi)\pi(x)),\quad x\in G.$$
If $f\in L^2(G)$ has the form
$f(x)=\tr(A_1\pi_1(x))+\cdots+\tr(A_k\pi_k(x))$, where
$[\pi_1],\ldots,[\pi_k]\in\hat{G}$ are distinct and
$A_i\in\MM(d_{\pi_i},\CC)$, then, by the Peter-Weyl Theorem, we have
$\supp(\hat{f})\subseteq\{[\pi_1],\ldots,[\pi_k]\}$ and
$\hat{f}(\pi_i)=A_i$. Here
$\supp(\hat{f})=\{[\pi]\in\hat{G}\mid\hat{f}(\pi)\neq0\}$.

Let $L^2_c(G)$ be the subspace of central functions in $L^2(G)$,
i.e., $f\in L^2_c(G)$ if and only if $f(xy)=f(yx)$ for almost all
$x,y\in G$. Then $f\in L^2_c(G)$ if and only if $\hat{f}(\pi)$ is a
scalar matrix for every $[\pi]\in\hat{G}$. Let $L^2_c(G)^\bot$ be
the orthogonal complement of $L^2_c(G)$ in $L^2(G)$. By the Fourier
inversion formula, one can show that $f\in L^2_c(G)^\bot$ if and
only if $\tr(\hat{f}(\pi))=0$ for every $[\pi]\in\hat{G}$.

A crucial property of the Fourier transform is that it converts the
regular representations of $G$ into matrix multiplications. As
usual, the left and right regular representations of
$G$ in $L^2(G)$ are defined by
$$(L_yf)(x)=f(y^{-1}x), \quad (R_yf)(x)=f(xy),$$ respectively, where
$f\in L^2(G)$ and $x,y\in G$.  Then it is easy to show that
$$(L_yf)\hat{}\,(\pi)=\hat{f}(\pi)\pi(y)^{-1}, \quad
(R_yf)\hat{}\,(\pi)=\pi(y)\hat{f}(\pi).$$

\subsection{Representation theory}

For a positive integer $n$, let $I_n$ denote the $n\times n$
identity matrix, and if $n$ is even, let
$J_n=\begin{bmatrix}0&I_{n/2}\\-I_{n/2}&0\\\end{bmatrix}$. If $n$ is
clear from the context, we will simply denote $I=I_n$ and $J=J_n$.
Recall that $Sp(n)=\{x\in U(n)\mid xJx^tJ^t=I\}$ if $n$ is even,
where $A^t$ refers to the transpose of a matrix
$A$. We recall the following definitions.

\begin{definition}
Let $\pi:G\rightarrow U(n)$ be an irreducible representation.
\begin{itemize}
\item[(1)] $\pi$ is of \emph{complex type} if $[\bar{\pi}]\neq[\pi]$.
\item[(2)] $\pi$ is of \emph{real type} if there exists $x\in U(n)$
such that $x\pi(G)x^{-1}\subseteq O(n)$.
\item[(3)] $\pi$ is of \emph{quaternionic type} if $n$ is even and
there exists $x\in U(n)$ such that $x\pi(G)x^{-1}\subseteq Sp(n)$.
\end{itemize}
\end{definition}

What is really important for us is the equivalence classes of
representations. So if $\pi$ is of real (resp. quaternionic) type,
we will always assume that $\pi(G)\subseteq O(n)$ (resp.
$\pi(G)\subseteq Sp(n)$).

Let $\hat{G}_c$ (resp. $\hat{G}_r$, $\hat{G}_q$) denote the set of
(equivalence classes of) irreducible representations of $G$ of
complex (resp. real, quaternionic) type. Then we have the following
basic fact.

\begin{theorem}\label{T:types}
$\hat{G}$ is the disjoint union of $\hat{G}_c$, $\hat{G}_r$, and
$\hat{G}_q$.
\end{theorem}

\begin{proof}[Proof (sketched).]
For an irreducible representation $\pi:G\rightarrow U(n)$, we
consider the representation $\rho$ of $G$ in $\MM(n,\CC)$ defined by
$\rho(g)(A)=\pi(g)A\pi(g)^t$. Let $\MM(n,\CC)^G$ denote the space of matrices $A$ such
that $\rho(g)(A)=A$ for all $g\in G$. Then $[\bar{\pi}]=[\pi]$ if
and only if $\dim\MM(n,\CC)^G=1$. In this case, any nonzero matrix
in $\MM(n,\CC)^G$ is invertible. It is easy to see that $\MM(n,\CC)$
is decomposed as the $G$-invariant direct sum of the space of
symmetric matrices $\MM_{\mathrm{symm}}(n,\CC)$ and the space of
skew-symmetric matrices $\MM_{\mathrm{skew}}(n,\CC)$. Hence
$[\bar{\pi}]=[\pi]$ if and only if either
$\dim\MM_{\mathrm{symm}}(n,\CC)^G=1$ (which means that $\pi(G)$ lies
in a conjugate of $O(n)$), or $\dim\MM_{\mathrm{skew}}(n,\CC)^G=1$
(which means that $n$ is even and $\pi(G)$ lies in a conjugate of
$Sp(n)$). Since $\dim\MM(n,\CC)^G=1$, the two cases can not occur
simultaneously. For more details, see \cite[Section 2.6]{BtD}.
\end{proof}

We define an equivalence relation on $\hat{G}$ for which the
equivalence class of $[\pi]$ is $\{[\pi],[\bar{\pi}]\}$ if
$[\pi]\in\hat{G}_c$, and $\{[\pi]\}$ if $[\pi]\in\hat{G}_r$ or
$\hat{G}_q$. We denote the equivalence class of $[\pi]$ with respect
to this equivalence relation by $[[\pi]]$, and the set of all
equivalence classes by $[\hat{G}]$.

\section{Constructing solutions from admissible tuples}

We first introduce some notions on solutions of
Eq.~\eqref{E:eq-general}, and examine their basic properties. For
$g,h\in L^2(G)$, let $g\otimes h$ be the function on $G^2$ defined
by $g\otimes h(x,y)=g(x)h(y)$. As being a solution
of Eq.~\eqref{E:eq-general} is a property about $f_1,f_2,f_3,f_4$
and $f_5\otimes f_6$, it is natural to denote a solution as a
$5$-tuple $\FF=(f_1,f_2,f_3,f_4,f_5\otimes f_6)$ of functions. But
sometimes we will also write the $5$-tuple $\FF$ as $(f_i)_{i=1}^6$
or simply $(f_i)$ for convenience.

The corresponding homogeneous equation \eqref{E:eq-homogeneous} is
important for us. Its solutions are $4$-tuples of functions
$(f_1,f_2,f_3,f_4)$, and form a closed subspace of $L^2(G)^4$ in the
usual way. If $(f_i)_{i=1}^6$ is a solution of
Eq.~\eqref{E:eq-general} satisfying $f_5\otimes f_6\equiv0$, then
$(f_i)_{i=1}^4$ is a solution of Eq.~\eqref{E:eq-homogeneous}. In
this case, without loss of generality, we always
assume that $f_5\equiv f_6\equiv0$. Conversely, if $(f_i)_{i=1}^4$
is a solution of Eq.~\eqref{E:eq-homogeneous}, then
$(f_1,f_2,f_3,f_4,0)$ is a solution of Eq.~\eqref{E:eq-general},
where $0$ is the zero function on $G^2$. We identify $(f_i)_{i=1}^4$
with $(f_1,f_2,f_3,f_4,0)$, and call such a solution a
\emph{homogeneous solution} of Eq.~\eqref{E:eq-general}. We say that
it is the \emph{trivial solution} if furthermore $f_i\equiv0$ for
$1\leq i\leq4$. If $\FF=(f_i)_{i=1}^6$ is a solution and
$\FF'=(f'_i)_{i=1}^4$ is a homogeneous solution, then their
\emph{sum} $\FF+\FF'=(f_1+f'_1,f_2+f'_2,f_3+f'_3,f_4+f'_4,f_5\otimes
f_6)$ is also a solution. It is obvious that if $c_1,c_2\in
L_c^2(G)$, then
\begin{equation}\label{E:c1c2}
\FF_{c_1,c_2}=(c_1,-c_1,c_2,-c_2)
\end{equation}
is a homogeneous solution. We say that a solution $(f_i)_{i=1}^6$ of
Eq.~\eqref{E:eq-general} is \emph{normalized} if $f_1-f_2,
f_3-f_4\in L_c^2(G)^\bot$. Then any solution of
Eq.~\eqref{E:eq-general} can be uniquely decomposed
as a sum $\FF+\FF_{c_1,c_2}$, where $\FF$ is normalized and
$\FF_{c_1,c_2}$ is given by \eqref{E:c1c2}. Furthermore, in the
Hilbert space of homogeneous solutions of Eq.~\eqref{E:eq-general},
normalized homogeneous solutions form the orthogonal complement of
the space of solutions of the form $\FF_{c_1,c_2}$. Finally, we say
that a solution $\FF=(f_i)_{i=1}^6$ of Eq.~\eqref{E:eq-general} is
\emph{pure} if $\bigcup_{i=1}^6\supp(\hat{f_i})\subseteq\varpi$ for
some $\varpi\in[\hat{G}]$. In this case, we say that $\FF$ is
\emph{supported on $\varpi$}.

In Section 5, we will determine all pure normalized solutions of
Eq.~\eqref{E:eq-general}, prove that pure
normalized homogeneous solutions span the space of normalized
homogeneous solutions, and that any solution is the sum of a pure
normalized solution and a homogeneous solution.

We will convert Eq.~\eqref{E:eq-general} into a family of matrix
equations. We call solutions of these matrix equations admissible
matrix tuples, whose definitions are as follows. For
$A,B,C,D,E,F\in\MM(n,\CC)$, we consider the linear maps
$\Phi_{A,B}^c$, $\Phi_{A,B,C,D}^r$, $\Phi_{A,B,C,D}^q$ (if $n$ is
even), and $\Psi_{E\otimes F}$ from $\MM(n,\CC)$ into itself defined
by
\begin{align*}
\Phi_{A,B}^c(X)&=AX+XB,\\
\Phi_{A,B,C,D}^r(X)&=AX+XB+(CX+XD)^t,\\
\Phi_{A,B,C,D}^q(X)&=AX+XB+J(CX+XD)^tJ^t,\\
\Psi_{E\otimes F}(X)&=\tr(EX)F.
\end{align*}
It is easy to see that $\Psi_{E\otimes F}$ depends only on $E\otimes
F\in\MM(n,\CC)\otimes\MM(n,\CC)$, and that if $n$ is even we have
\begin{equation}\label{E:FIJ}
\Phi^r_{A,B,C,D}(X)=-\Phi^q_{A,-JBJ,-C,JDJ}(XJ)J.
\end{equation}

\begin{definition}\label{D:adm}
Let $A,B,C,D,E,F\in\MM(n,\CC)$.
\begin{itemize}
\item[(1)] $(A,B,E\otimes F)$ is an \emph{admissible tuple of complex type}
(\emph{$c$-admissible tuple} abbreviated) if $\tr(A)=\tr(B)$ and
$\Phi_{A,B}^c=\Psi_{E\otimes F}$.
\item[(2)] $(A,B,C,D,E\otimes F)$ is an \emph{admissible tuple of real
type} (\emph{$r$-admissible tuple} abbreviated) if $\tr(A)=\tr(B)$,
$\tr(C)=\tr(D)$, and $\Phi_{A,B,C,D}^r=\Psi_{E\otimes F}$.
\item[(3)] $(A,B,C,D,E\otimes F)$ is an \emph{admissible tuple of
quaternionic type} (\emph{$q$-admissible tuple} abbreviated) if $n$
is even, $\tr(A)=\tr(B)$, $\tr(C)=\tr(D)$, and
$\Phi_{A,B,C,D}^q=\Psi_{E\otimes F}$.
\end{itemize}
\end{definition}

We refer to $n$ as the \emph{order} of the above admissible matrix
tuples. An admissible tuple $\TT$ is \emph{homogeneous} if $E\otimes
F=0$, and is \emph{trivial} if $A=B(=C=D)=0$. It is obvious that
trivial admissible tuples are homogeneous. If $\TT$ is homogeneous,
we always assume that $E=F=0$.

We should mention that the trace
conditions in Definition \ref{D:adm} are not essential. As we will
see later, they are imposed so that admissible tuples correspond to
normalized solutions. This will simplify some arguments below.

We will determine all admissible matrix tuples in the next section.
In the rest of this section, we explain how to construct pure
normalized solutions of Eq.~\eqref{E:eq-general} from admissible
tuples. We also exhibit some examples of admissible tuples, which
indeed include all nontrivial ones. The solutions
constructed from these examples form the building blocks of the
general solution of Eq.~\eqref{E:eq-general}.

We begin with a simple example.

\begin{example}\label{Ex:U1}
Let $\varepsilon_1,\delta_1,\varepsilon_2,\delta_2\in\CC$. Then
$(\varepsilon_i\delta_j/2,\varepsilon_i\delta_j/2,\varepsilon_i
\delta_j)$ $(i,j=1,2)$ are $1$-ordered $c$-admissible tuples. Define
the tuple of functions
$\FF_{\varepsilon_1,\delta_1,\varepsilon_2,\delta_2}^{U(1)}=(f_i)_{i=1}^6$
as
$$\FF_{\varepsilon_1,\delta_1,\varepsilon_2,\delta_2}^{U(1)}: \quad \begin{cases}
f_1(x)=f_2(x)=(\varepsilon_1\delta_1x+\varepsilon_2\delta_2\bar{x})/2,\\
f_3(x)=f_4(x)=(\varepsilon_1\delta_2x+\varepsilon_2\delta_1\bar{x})/2,\\
f_5\otimes
f_6(x,y)=(\varepsilon_1x+\varepsilon_2\bar{x})(\delta_1y+\delta_2\bar{y}),
\end{cases} \quad x,y\in U(1).$$
Then it is easy to check that
$\FF_{\varepsilon_1,\delta_1,\varepsilon_2,\delta_2}^{U(1)}$ is a
pure normalized solution of Eq.~\eqref{E:eq-general} on $U(1)$
supported on $[[\iota_{U(1)}]]$, where $\iota_{U(1)}$ is the
identity representation of $U(1)$. It is homogeneous if and only if
it is the trivial solution.
\end{example}

The general principle of constructing solutions from admissible
tuples of real and quaternionic types is as follows. For a closed
irreducible subgroup $K$ of $U(n)$ and a matrix $L\in\MM(n,\CC)$, we
define the function $f_L$ on $K$ as $f_L(x)=\tr(Lx)$, $x\in K$. Then
we have $\supp(\hat{f}_L)\subseteq\{[\iota_K]\}$ and
$\hat{f}_L(\iota_K)=L$, where $\iota_K:K\rightarrow U(n)$ is the
inclusion. For a $5$-tuple $\TT=(A,B,C,D,E\otimes F)$, where
$A,\ldots,F\in\MM(n,\CC)$, we define the $5$-tuple of functions
$$\FF_\TT^K=(f_A,f_B,f_C,f_D,f_E\otimes f_F).$$ Clearly, $f_E\otimes
f_F$ depends only on $E\otimes F$.

\begin{proposition}\label{P:construction}
We keep the notation as above.
\begin{itemize}
\item[(1)] If $\TT$ is an $n$-ordered $r$-admissible tuple, then
$\FF_\TT^{O(n)}$ is a pure normalized solution of
Eq.~\eqref{E:eq-general} on $O(n)$ supported on
$\{[\iota_{O(n)}]\}$. $\FF_\TT^{O(n)}$ is homogeneous if and only if
$\TT$ is homogeneous.
\item[(2)] If $n$ is even and $\TT$ is an $n$-ordered $q$-admissible
tuple, then $\FF_\TT^{Sp(n)}$ is a pure normalized solution of
Eq.~\eqref{E:eq-general} on $Sp(n)$ supported on
$\{[\iota_{Sp(n)}]\}$. $\FF_\TT^{Sp(n)}$ is homogeneous if and only
if $\TT$ is homogeneous.
\end{itemize}
\end{proposition}

\begin{proof}
(1) Since $\Phi^r_{A,B,C,D}=\Psi_{E\otimes F}$, for all
$x,y\in O(n)$ we have
\begin{align*}
&f_A(xy)+f_B(yx)+f_C(xy^{-1})+f_D(y^{-1}x)\\
=&\tr(Axy)+\tr(Byx)+\tr(Cxy^t)+\tr(Dy^tx)\\
=&\tr(Axy+xBy+x^tC^ty+D^tx^ty)\\
=&\tr(\Phi^r_{A,B,C,D}(x)y)\\
=&\tr(\tr(Ex)Fy)\\
=&f_E(x)f_F(y).
\end{align*}
So $\FF_\TT^{O(n)}$ is a solution of Eq.~\eqref{E:eq-general} on
$O(n)$. Obviously it is a pure solution supported on $\{[\iota]\}$,
where $\iota=\iota_{O(n)}$. Since
$\tr(\hat{f}_A(\iota)-\hat{f}_B(\iota))=\tr(A-B)=0$, we have
$f_A-f_B\in L_c^2(O(n))^\bot$. Similarly, $f_C-f_D\in
L_c^2(O(n))^\bot$. Thus $\FF_\TT^{O(n)}$ is normalized. It is
homogeneous if and only if $f_E\equiv0$ or $f_F\equiv0$, which is
equivalent to $E\otimes
F=\hat{f}_E(\iota)\otimes\hat{f}_F(\iota)=0$, i.e., $\TT$ is
homogeneous.

(2) Since $\Phi^q_{A,B,C,D}=\Psi_{E\otimes F}$, for all $x,y\in
Sp(n)$ we have
\begin{align*}
&f_A(xy)+f_B(yx)+f_C(xy^{-1})+f_D(y^{-1}x)\\
=&\tr(Axy)+\tr(Byx)+\tr(CxJy^tJ^t)+\tr(DJy^tJ^tx)\\
=&\tr(Axy+xBy+J^tx^tC^tJy+J^tD^tx^tJy)\\
=&\tr(\Phi^q_{A,B,C,D}(x)y)\\
=&\tr(\tr(Ex)Fy)\\
=&f_E(x)f_F(y).
\end{align*}
Hence $\FF_\TT^{Sp(n)}$ is a solution of Eq.~\eqref{E:eq-general} on
$Sp(n)$. The proofs of the other assertions in (2) are similar to
those of the corresponding parts in (1) and omitted here.
\end{proof}

Note that if $\varphi:G\rightarrow K$ is a homomorphism and
$\FF^K=(f_i)$ is a solution of Eq.~\eqref{E:eq-general} on $K$, then
$\FF^K\circ\varphi=(f_i\circ\varphi)$ is a solution on $G$. Some
relations between $\FF^K$ and $\FF^K\circ\varphi$ are revealed in
the following assertion.

\begin{proposition}\label{P:construction2}
Let $\pi:G\rightarrow U(n)$ be an irreducible representation of
complex (resp. real, quaternionic) type, and let $K=U(n)$ (resp.
$O(n)$, $Sp(n)$). If $\FF^K=(f_i)$ is a pure solution of
Eq.~\eqref{E:eq-general} on $K$ supported on $[[\iota_K]]$, then
$\FF^K\circ\pi$ is a pure solution of Eq.~\eqref{E:eq-general} on
$G$ supported on $[[\pi]]$, and $\FF^K\circ\pi$ is normalized (resp.
homogeneous) if and only if $\FF^K$ is normalized (resp.
homogeneous).
\end{proposition}

\begin{proof}
It suffices to prove that if $f$ is a function on $K$ with
$\supp(\hat{f})\subseteq[[\iota_K]]$, then
$\supp((f\circ\pi)\hat{}\,)\subseteq[[\pi]]$, $f\circ\pi\in
L^2_c(G)^\bot$ if and only if $f\in L^2_c(K)^\bot$, and
$f\circ\pi\equiv0$ if and only if $f\equiv0$. Suppose that $\pi$ is
of complex type. Then $f$ is of the form
$f(x)=\tr(Ax)+\tr(B\bar{x})$, where $x\in U(n)$, $A,B\in\MM(n,\CC)$.
Hence $(f\circ\pi)(y)=\tr(A\pi(y))+\tr(B\bar{\pi}(y))$, $y\in G$.
This implies that $(f\circ\pi)\hat{}\,(\pi)=A$,
$(f\circ\pi)\hat{}\,(\bar{\pi})=B$, and
$(f\circ\pi)\hat{}\,(\pi')=0$ if $[\pi']\notin[[\pi]]$. So
$\supp((f\circ\pi)\hat{}\,)\subseteq[[\pi]]$. Moreover, we have
$$f\circ\pi\in L^2_c(G)^\bot\Leftrightarrow\tr A=\tr B=0\Leftrightarrow f\in L^2_c(K)^\bot,$$
$$f\circ\pi\equiv0\Leftrightarrow A=B=0\Leftrightarrow f\equiv0.$$ The proofs of the other
two cases are similar and left to the reader.
\end{proof}

\begin{example}\label{Ex:O1}
Any $1$-ordered $r$-admissible tuple is of the form
$$\TT_{a,b}=(a/2,a/2,b/2,b/2,a+b)$$ for some $a,b\in\CC$. It is
homogeneous if and only if $a+b=0$. We define the
tuple of functions $\FF_{a,b}^{O(1)}=(f_i)_{i=1}^6$ as
$$\FF_{a,b}^{O(1)}: \quad \begin{cases}
f_1(x)=f_2(x)=ax/2, \\
f_3(x)=f_4(x)=bx/2,\\
f_5\otimes f_6(x,y)=(a+b)xy,
\end{cases} \quad x,y\in O(1).$$
Then $\FF_{a,b}^{O(1)}=\FF^{O(1)}_{\TT_{a,b}}$. By Proposition
\ref{P:construction} (1), it is a pure normalized solution of
Eq.~\eqref{E:eq-general} on $O(1)$ supported on
$\{[\iota_{O(1)}]\}$. It is homogeneous if and only if $a+b=0$.
Note that $\FF_{a,b}^{O(1)}$ is the restriction of the solution
$\FF_{a,1,b,0}^{U(1)}$ on $U(1)$ (see Example \ref{Ex:U1}). But it
may occur that $\FF_{a,1,b,0}^{U(1)}$ is non-homogeneous while
$\FF_{a,b}^{O(1)}$ is homogeneous. This fact is meaningful when we
construct the general solution of Eq.~\eqref{E:eq-general} on
arbitrary compact groups (see Section 5). For later reference, we
denote $\FF_a^{O(1)}=\FF_{2a,-2a}^{O(1)}$. In our notation of
homogeneous solutions, $\FF_a^{O(1)}$ is the $4$-tuple of functions
$(f_i)_{i=1}^4$ defined as
$$\FF_a^{O(1)}: \quad f_1(x)=f_2(x)=-f_3(x)=-f_4(x)=ax, \quad x\in O(1).$$
\end{example}

Now we consider admissible matrix tuples of higher order. Since the
bilinear pairing $(X,Y)\mapsto\tr(XY)$ on $\MM(n,\CC)$ is
non-degenerate, for a linear map
$\Gamma:\MM(n,\CC)\rightarrow\MM(n,\CC)$, we can define its
\emph{adjoint} $\Gamma^\dag$ by
$\tr(\Gamma(X)Y)=\tr(X\Gamma^\dag(Y))$ for all $X,Y\in\MM(n,\CC)$.
It is straightforward to check that
\begin{align}
(\Phi_{A,B}^c)^\dag&=\Phi_{B,A}^c,\label{E:adj-Phi-c}\\
(\Phi_{A,B,C,D}^r)^\dag&=\Phi_{B,A,C^t,D^t}^r,\label{E:adj-Phi-r}\\
(\Phi_{A,B,C,D}^q)^\dag&=\Phi_{B,A,JC^tJ^t,JD^tJ^t}^q, \label{E:adj-Phi-q}\\
(\Psi_{E\otimes F})^\dag&=\Psi_{F\otimes E}.\label{E:adj-Psi}
\end{align}

\begin{lemma}\label{L:n=2}
Let $A,B\in\MM(2,\CC)$ be such that $\tr(A)=\tr(B)$.
\begin{itemize}
\item[(1)] The tuples
\begin{align*}
\TT^r_{A,B}&=(A,B,-A,-B,-(JA+BJ)\otimes J),\\
(\TT^r_{A,B})^\dag&=(A,B,-B^t,-A^t,-J\otimes(AJ+JB))
\end{align*}
are $r$-admissible. They are homogeneous if and only if $\tr(A)=0$
and $B=A^t$.
\item[(2)] The tuples
\begin{align*}
\TT^q_{A,B}&=(A,B,A,B,(A+B)\otimes I),\\
(\TT^q_{A,B})^\dag&=(A,B,JB^tJ^t,JA^tJ^t,I\otimes(A+B))
\end{align*}
are $q$-admissible. They are homogeneous if and only if $\tr(A)=0$
and $B=-A$.
\end{itemize}
\end{lemma}

\begin{proof}
We first prove the assertions for $\TT^r_{A,B}$ and
$\TT^q_{A,B}$. Since $Y+JY^tJ^t=\tr(Y)I$ for any $Y\in\MM(2,\CC)$,
we have
\begin{align}\label{E:q}
\Phi^q_{A,B,A,B}(X)&=AX+XB+J(AX+XB)^tJ^t \notag\\
&=\tr(AX+XB)I\notag\\
&=\Psi_{(A+B)\otimes I}(X).
\end{align}
So $\TT^q_{A,B}$ is $q$-admissible. By \eqref{E:FIJ}, we have
\begin{align*}
\Phi^r_{A,B,-A,-B}(X)&=-\Phi^q_{A,-JBJ,A,-JBJ}(XJ)J\\
&=-\Psi_{(A-JBJ)\otimes I}(XJ)J\\
&=\Psi_{-(JA+BJ)\otimes J}(X).
\end{align*}
So $\TT^r_{A,B}$ is $r$-admissible.

Now by
\eqref{E:adj-Phi-r}--\eqref{E:adj-Psi}, we have
\begin{equation}\label{E:qdag}
\Phi^q_{A,B,JB^tJ^t,JA^tJ^t}=(\Phi^q_{B,A,B,A})^\dag=(\Psi_{(A+B)\otimes
I})^\dag=\Psi_{I\otimes(A+B)},
\end{equation}
$$\Phi^r_{A,B,-B^t,-A^t}=(\Phi^r_{B,A,-B,-A})^\dag=(\Psi_{-(JB+AJ)\otimes J})^\dag=\Psi_{-J\otimes(AJ+JB)}.$$
Hence $(\TT^q_{A,B})^\dag$ and $(\TT^r_{A,B})^\dag$
are admissible tuples of quaternionic and real type, respectively.

The conditions of being homogeneous are easy to prove and left to
the reader.
\end{proof}

\begin{remark}\label{R:3.1}
The families $\TT^r_{A,B}$ and $(\TT^r_{A,B})^\dag$ (resp.
$\TT^q_{A,B}$ and $(\TT^q_{A,B})^\dag$) are not mutually exclusive.
Indeed, it is easy to check that $\TT^r_{A,B}=(\TT^r_{A,B})^\dag$ if
and only if $B=A^t$,  and $\TT^q_{A,B}=(\TT^q_{A,B})^\dag$ if and
only if $B=\tr(A)I-A$. In particular, if $\TT^r_{A,B}$ (or,
equivalently, $(\TT^r_{A,B})^\dag$) is homogeneous, then
$\TT^r_{A,B}=(\TT^r_{A,B})^\dag$. Similarly, if $\TT^q_{A,B}$ (or
$(\TT^q_{A,B})^\dag$) is homogeneous, then
$\TT^q_{A,B}=(\TT^q_{A,B})^\dag$.
\end{remark}

\begin{example}\label{Ex:O2}
Let $A,B\in\MM(2,\CC)$ with $\tr(A)=\tr(B)$. By Proposition
\ref{P:construction} (1) and Lemma \ref{L:n=2} (1), the tuples of
functions $\FF^{O(2)}_{A,B}=\FF^{O(2)}_{\TT^r_{A,B}}$ and
$(\FF^{O(2)}_{A,B})^\dag=\FF^{O(2)}_{(\TT^r_{A,B})^\dag}$ are pure
normalized solutions of Eq.~\eqref{E:eq-general} on $O(2)$ supported
on $\{[\iota_{O(2)}]\}$. Writing explicitly, we have
$$\FF^{O(2)}_{A,B}: \quad \begin{cases}
f_1(x)=-f_3(x)=\tr(Ax),\\
f_2(x)=-f_4(x)=\tr(Bx),\\
f_5\otimes f_6(x,y)=-\tr((JA+BJ)x)\tr(y),
\end{cases} \quad x,y\in O(2);$$
$$(\FF^{O(2)}_{A,B})^\dag: \quad \begin{cases}
f_1(x)=-f_4(x^{-1})=\tr(Ax),\\
f_2(x)=-f_3(x^{-1})=\tr(Bx),\\
f_5\otimes f_6(x,y)=-\tr(x)\tr((AJ+JB)y),
\end{cases} \quad x,y\in O(2).$$
The solutions $\FF^{O(2)}_{A,B}$ and $(\FF^{O(2)}_{A,B})^\dag$ are
homogeneous if and only if $\tr(A)=0$ and $B=A^t$. In this case the
two solutions are equal (see Remark \ref{R:3.1}). We denote
$\FF^{O(2)}_A=\FF^{O(2)}_{A,A^t}=(\FF^{O(2)}_{A,A^t})^\dag$ if
$\tr(A)=0$. The functions in $\FF^{O(2)}_A$ are
$$\FF^{O(2)}_A: \quad f_1(x)=f_2(x^{-1})=-f_3(x)=-f_4(x^{-1})=\tr(Ax), \quad x\in O(2).$$
\end{example}

\begin{example}\label{Ex:SU2}
Let $A,B\in\MM(2,\CC)$ with $\tr(A)=\tr(B)$. By Proposition
\ref{P:construction} (2), Lemma \ref{L:n=2} (2) and the fact
$Sp(2)=SU(2)$, the tuples of functions
$\FF^{SU(2)}_{A,B}=\FF^{SU(2)}_{\TT^q_{A,B}}$ and
$(\FF^{SU(2)}_{A,B})^\dag=\FF^{SU(2)}_{(\TT^q_{A,B})^\dag}$ are pure
normalized solutions of Eq.~\eqref{E:eq-general} on $SU(2)$
supported on $\{[\iota_{SU(2)}]\}$. The functions in these
solutions are
$$
\FF^{SU(2)}_{A,B}: \quad \begin{cases}
f_1(x)=f_3(x)=\tr(Ax),\\
f_2(x)=f_4(x)=\tr(Bx),\\
f_5\otimes f_6(x,y)=\tr((A+B)x)\tr(y),
\end{cases} \quad x,y\in SU(2);
$$
$$
(\FF^{SU(2)}_{A,B})^\dag: \quad \begin{cases}
f_1(x)=\tr(Ax),\\
f_2(x)=\tr(Bx),\\
f_3(x)=\tr(A)\tr(x)-f_2(x),\\
f_4(x)=\tr(A)\tr(x)-f_1(x),\\
f_5\otimes f_6(x,y)=\tr(x)\tr((A+B)y),
\end{cases} \quad x,y\in SU(2).
$$
These solutions are homogeneous if and only if $\tr(A)=0$ and
$B=-A$, and in this case we have
$\FF^{SU(2)}_{A,B}=(\FF^{SU(2)}_{A,B})^\dag$. We denote
$\FF^{SU(2)}_A=\FF^{SU(2)}_{A,A^t}=(\FF^{SU(2)}_{A,A^t})^\dag$ if
$\tr(A)=0$. Writing explicitly, it is
$$\FF^{SU(2)}_A: \quad f_1(x)=-f_2(x)=f_3(x)=-f_4(x)=\tr(Ax), \quad x\in SU(2).$$
\end{example}

Now we consider $3$-ordered $r$-admissible tuples. We view elements
of $\CC^3$ as column vectors. For $u,v\in\CC^3$, let $\langle
u,v\rangle=u^tv$ be the standard bilinear pairing, and define
$$
\tau_{u,v}=uv^t-\frac{1}{2}\langle u,v\rangle I_3\in\MM(3,\CC).
$$
Let $\MM_{\mathrm{skew}}(3,\CC)$ denote the space of $3\times3$
skew-symmetric complex matrices. For
$u=(u_1,u_2,u_3)^t\in\CC^3$, let
$$
\sigma_u=\begin{bmatrix}0&-u_3&u_2\\u_3&0&-u_1\\-u_2&u_1&0\end{bmatrix}\in\MM_{\mathrm{skew}}(3,\CC).
$$
Note that for $w\in\CC^3$, $\sigma_uw$ is (the complex analogue of)
the cross product $u\times w$ of $u$ and $w$.

\begin{lemma}\label{L:n=3}
For any $u,v\in\CC^3$, the tuple
$$\TT_{u,v}=(\tau_{u,v},\tau_{v,u},-\tau_{u,v},-\tau_{v,u},\sigma_u\otimes\sigma_v)$$
is $r$-admissible. It is homogeneous if and only if it is the
trivial tuple.
\end{lemma}

\begin{proof}
Firstly we consider the representations $\rho_1$ and $\rho_2$ of the
Lie algebra $\Lgl(3,\CC)$ in $\MM_{\mathrm{skew}}(3,\CC)$ and
$\CC^3$ defined by
\begin{equation}\label{E:representation}
\rho_1(A)(Y)=AY+YA^t, \quad \rho_2(A)(w)=(\tr(A)I_3-A^t)w,
\end{equation}
respectively, where $A\in\Lgl(3,\CC)$,
$Y\in\MM_{\mathrm{skew}}(3,\CC)$, $w\in\CC^3$. We claim that the
linear isomorphism
$\sigma:\CC^3\rightarrow\MM_{\mathrm{skew}}(3,\CC)$ sending $w$ to
$\sigma_w$ is an equivalence between $\rho_1$ and $\rho_2$, i.e.,
\begin{equation}\label{E:equivalence}
\rho_1(A)(\sigma_w)=\sigma(\rho_2(A)(w))
\end{equation}
for all $A\in\Lgl(3,\CC)$ and $w\in\CC^3$. To prove this, we note
(the complex analogue of) the equality for scalar triple products,
i.e., for all $w,w_1,w_2\in\CC^3$, we have
$$\langle\sigma_ww_1,w_2\rangle=\det[w,w_1,w_2],$$
where $[w,w_1,w_2]$ is the $3\times3$ matrix specified by column
vectors. Now let $A\in\Lgl(3,\CC)$ and $w,w_1,w_2\in\CC^3$. Then we
have
\begin{align*}
\langle\rho_1(A)(\sigma_w)w_1,w_2\rangle
&=\langle (A\sigma_w+\sigma_wA^t)w_1,w_2\rangle\\
&=\langle A\sigma_ww_1,w_2\rangle + \langle\sigma_wA^tw_1,w_2\rangle\\
&=\langle \sigma_ww_1,A^tw_2\rangle + \langle\sigma_wA^tw_1,w_2\rangle\\
&=\det[w,w_1,A^tw_2]+\det[w,A^tw_1,w_2]
\end{align*}
and
\begin{align*}
\langle\sigma(\rho_2(A)(w))w_1,w_2\rangle
&=\det[\rho_2(A)(w),w_1,w_2]\\
&=\det[(\tr(A)I_3-A^t)w,w_1,w_2]\\
&= \tr(A)\det[w,w_1,w_2]-\det[A^tw,w_1,w_2].
\end{align*}
This proves \eqref{E:equivalence} by noting the fact that
$$\det[Aw,w_1,w_2]+\det[w,Aw_1,w_2]+\det[w,w_1,Aw_2]=\tr(A)\det[w,w_1,w_2]$$
for all $A\in\Lgl(3,\CC)$ and $w,w_1,w_2\in\CC^3$.

Now we notice that
\begin{align*}
\rho_2(\tau_{u,v})(w)&=-\langle u,w\rangle
v=\frac{1}{2}\tr(\sigma_u\sigma_w)v, \\
\tau_{u,v}^t&=\tau_{v,u}, \quad \sigma_u^t=-\sigma_u.
\end{align*}
From these identities, \eqref{E:representation}, and \eqref{E:equivalence},
it follows that for all $X\in\MM(3,\CC)$
we have
\begin{align}\label{E:tausigma}
&\Phi_{\tau_{u,v},\tau_{v,u},-\tau_{u,v},-\tau_{v,u}}^r(X)\notag\\
=&\tau_{u,v}X+X\tau_{v,u}-(\tau_{u,v}X+X\tau_{v,u})^t\notag\\
=&\tau_{u,v}(X-X^t)+(X-X^t)\tau_{u,v}^t\notag\\
=&\rho_1(\tau_{u,v})(X-X^t)=\sigma(\rho_2(\tau_{u,v})(\sigma^{-1}(X-X^t)))\notag\\
=&-\sigma(\langle u,\sigma^{-1}(X-X^t)\rangle v)=-\langle
u,\sigma^{-1}(X-X^t)\rangle\sigma(v)\notag\\
=&\frac{1}{2}\tr(\sigma_u(X-X^t))\sigma_v=\tr(\sigma_uX)\sigma_v\notag\\
=&\Psi_{\sigma_u\otimes\sigma_v}(X).
\end{align}
This proves that $\TT_{u,v}$ is $r$-admissible. If $\TT_{u,v}$ is
homogeneous, then $\sigma_u=0$ or $\sigma_v=0$, which implies that
$u=0$ or $v=0$. Hence it is the trivial tuple.
\end{proof}

\begin{example}\label{Ex:O3}
For $u,v\in\CC^3$, we define the tuple of functions
$\FF_{u,v}^{O(3)}$ as $\FF^{O(3)}_{\TT_{u,v}}$. The functions in
$\FF_{u,v}^{O(3)}$ are
$$\FF_{u,v}^{O(3)}: \quad \begin{cases}
f_1(x)=f_2(x^{-1})=-f_3(x)=-f_4(x^{-1})=\tr(\tau_{u,v}x),\\
f_5\otimes f_6(x,y)=\tr(\sigma_ux)\tr(\sigma_vy),
\end{cases} \quad x,y\in O(3).$$
Then Proposition
\ref{P:construction} (1) and Lemma \ref{L:n=3} imply that
$\FF_{u,v}^{O(3)}$ is a pure normalized solution of
Eq.~\eqref{E:eq-general} on $O(3)$ supported on
$\{[\iota_{O(3)}]\}$. It is homogeneous if and only if it is the
trivial solution.
\end{example}

\section{Determination of admissible tuples}

In this section we determine all admissible matrix tuples, which are
completely described in the following three propositions. We keep
the same notation from Section 3.

\begin{proposition}\label{P:c}
Let $\TT=(A,B,E\otimes F)$ be an $n$-ordered $c$-admissible tuple.
\begin{itemize}
\item[(1)] If $n=1$, then $\TT=(a,a,2a)$ for some $a\in\CC$.
\item[(2)] If $n\geq2$, then $\TT$ is the trivial tuple.
\end{itemize}
\end{proposition}

\begin{proposition}\label{P:r}
Let $\TT=(A,B,C,D,E\otimes F)$ be an $n$-ordered $r$-admissible
tuple.
\begin{itemize}
\item[(1)] If $n=1$, then $\TT=\TT_{a,b}$ for some $a,b\in\CC$.
\item[(2)] If $n=2$, then $\TT=\TT^r_{A,B}$ or $(\TT^r_{A,B})^\dag$ for
some $A,B\in\MM(2,\CC)$ with $\tr(A)=\tr(B)$.
\item[(3)] If $n=3$, then $\TT=\TT_{u,v}$ for some $u,v\in\CC^3$.
\item[(4)] If $n\geq4$, then $\TT$ is the trivial tuple.
\end{itemize}
\end{proposition}

\begin{proposition}\label{P:q}
Let $n$ be even, and let $\TT=(A,B,C,D,E\otimes F)$ be an
$n$-ordered $q$-admissible tuple.
\begin{itemize}
\item[(1)] If $n=2$, then $\TT=\TT^q_{A,B}$ or $(\TT^q_{A,B})^\dag$ for
some $A,B\in\MM(2,\CC)$ with $\tr(A)=\tr(B)$.
\item[(2)] If $n\geq4$, then $\TT$ is the trivial tuple.
\end{itemize}
\end{proposition}

The assertions in Propositions \ref{P:c} (1) and \ref{P:r} (1) are
trivial. It remains to prove the others. Since our proofs of
\ref{P:q} (2) and \ref{P:r} (2) make use of \ref{P:r} (4) and
\ref{P:q} (1), respectively, and the proofs of \ref{P:c} (2) and
\ref{P:r} (4) are similar, we proceed the proofs in the following
order:
$$\text{\ref{P:c} (2), \ref{P:r} (4) $\Rightarrow$ \ref{P:q} (2), \ref{P:q} (1) $\Rightarrow$
\ref{P:r} (2), \ref{P:r} (3).}$$

\begin{proof}[Proof of Proposition \ref{P:c} (2)]
Denote $\Phi=\Phi_{A,B}^c$ and $\NN_n=\{1,\ldots,n\}$. Since
$\Phi=\Psi_{E\otimes F}$, we have $\dim\im(\Phi)\leq1$. So the
entries $\Phi(X)_{ij}$ $(i,j\in\NN_n)$ of $\Phi(X)$, viewed as
linear polynomials in the entries $X_{ij}$ of $X$, are mutually
linearly dependent. We make the convention that if a linear
polynomial $p$ in the variables $y_1,\ldots,y_m$ is written in the
reduced form as $p(y)=a_1y_1+a_2y_2+\cdots$, then the terms being
omitted do not contain $y_1$ and $y_2$.

Let $i,j\in\NN_n$, $i\neq j$. It is easy to see that
\begin{align*}
\Phi(X)_{ii}&=A_{ij}X_{ji}+0X_{jj}+\cdots,\\
\Phi(X)_{ij}&=0X_{ji}+A_{ij}X_{jj}+\cdots.
\end{align*}
Since they are linearly dependent, we must have $A_{ij}=0$. So $A$
is diagonal. Similarly, $B$ is diagonal. Now we have
$$\Phi(X)_{rs}=(A_{rr}+B_{ss})X_{rs}\quad \text{for all}\quad r,s\in\NN_n.$$
Setting $(r,s)=(i,i),(i,j),(j,i),(j,j)$, we get four polynomials.
Their mutual linear dependence implies that at most one of the four
sums $A_{ii}+B_{ii}$, $A_{ii}+B_{jj}$, $A_{jj}+B_{ii}$,
$A_{jj}+B_{jj}$ is nonzero. This forces that they are all zero. So
$A=-B\in\CC I$. But we have $\tr(A)=\tr(B)$. Hence $A=B=0$. This
proves that $\TT$ is the trivial tuple.
\end{proof}

We use the similar idea to prove \ref{P:r} (4).

\begin{proof}[Proof of Proposition \ref{P:r} (4)]
Denote $\Phi=\Phi_{A,B,C,D}^r$. Then $\dim\im(\Phi)\leq1$ and
$\Phi(X)_{ij}$ $(i,j\in\NN_n)$ are mutually linearly dependent. Let
$i,j\in\NN_n$ with $i\neq j$. Since $n\geq4$, there exist
$k,l\in\NN_n$ such that $i,j,k,l$ are distinct. Then we compute
\begin{align*}
\Phi(X)_{ik}&=A_{ij}X_{jk}+0X_{jl}+\cdots,\\
\Phi(X)_{il}&=0X_{jk}+A_{ij}X_{jl}+\cdots.
\end{align*}
Since they are linearly dependent, we have $A_{ij}=0$. So $A$ is
diagonal. Similarly, $B,C,D$ are diagonal. Now we have
$$\Phi(X)_{rs}=(A_{rr}+B_{ss})X_{rs}+(C_{ss}+D_{rr})X_{sr}\quad \text{for all}\quad
r,s\in\NN_n.$$ Setting $(r,s)=(i,j),(i,l),(k,j),(k,l)$, we get four
polynomials. Their mutual linear dependence implies
that at most one of $A_{ii}+B_{jj}$, $A_{ii}+B_{ll}$,
$A_{kk}+B_{jj}$, $A_{kk}+B_{ll}$ is nonzero. This forces that they
are all zero. So $A_{ii}+B_{jj}=0$ whenever $i\neq j$. This is
impossible unless $A=-B\in\CC I$. But we have $\tr(A)=\tr(B)$. So
$A=B=0$. Similarly, $C=D=0$. Hence $\TT$ is trivial.
\end{proof}

We now use Proposition \ref{P:r} (4) to prove Proposition \ref{P:q} (2).

\begin{proof}[Proof of Proposition \ref{P:q} (2)]
Suppose $n\geq4$ and $\TT$ is $q$-admissible. Then it follows from
\eqref{E:FIJ} that the tuple $(A,-JBJ,-C,JDJ,-(JE)\otimes(FJ))$ is
$r$-admissible. By Proposition \ref{P:r} (4), we have
$A=-JBJ=-C=JDJ=0$. So $A=B=C=D=0$.
\end{proof}

Similarly, due to \eqref{E:FIJ}, Proposition \ref{P:r} (2) is
equivalent to Proposition \ref{P:q} (1). We find that the proof of
Proposition \ref{P:q} (1) is easier to write up. So we prove it
first. In the following proof, we will constantly use the fact that
$Y+JY^tJ^t=\tr(Y)I$ for all $Y\in\MM(2,\CC)$
without any further mention.

\begin{proof}[Proof of Proposition \ref{P:q} (1)]
Denote $\Phi=\Phi_{A,B,C,D}^q$. Then $\dim\im(\Phi)\leq1$ and
$\Phi(X)_{ij}$ $(i,j\in\NN_2)$ are mutually linearly dependent. We
divide the proof into two steps.

\emph{Step (i).} First we assume that $C=-A$ and $D=-B$. We prove
that $\tr(A)=0$, $B=A$, and $\Phi(X)=2\tr(X)A$.

In this case, we have
$$\Phi(X)=AX+XB-J(AX+XB)^tJ^t.$$
Let $(i,j)=(1,2)$ or $(2,1)$. Since
\begin{align*}
\Phi(X)_{ii}&=(A_{ij}-B_{ij})X_{ji}+\cdots,\\
\Phi(X)_{ij}&=0X_{ji}+2A_{ij}X_{jj}+2B_{ij}X_{ii}+\cdots,
\end{align*}
their linear dependence implies that $A_{ij}=B_{ij}$. Using this, it
is easy to compute that
\begin{align*}
\Phi(X)_{ij}&=2(A_{ii}+B_{jj})X_{ij}+\cdots,\\
\Phi(X)_{ji}&=0X_{ij}+2(A_{jj}+B_{ii})X_{ji}+\cdots,\\
\Phi(X)_{ii}&=(A_{ii}+B_{ii})X_{ii}-(A_{jj}+B_{jj})X_{jj}.
\end{align*}
We claim that $A_{ii}+B_{jj}=0$. For otherwise, if
$A_{ii}+B_{jj}\neq0$, then by the mutual linear dependence, we have
$A_{jj}+B_{ii}=A_{ii}+B_{ii}=A_{jj}+B_{jj}=0$, which conflicts with
$A_{ii}+B_{jj}\neq0$. Now if
$A=\begin{bmatrix}a&b\\c&d\\\end{bmatrix}$, then
$B=\begin{bmatrix}-d&b\\c&-a\\\end{bmatrix}$. But $\tr(A)=\tr(B)$.
Hence $\tr(A)=0$ and $B=A$. This also implies that $C=D=-A=JA^tJ^t$.
By \eqref{E:qdag}, we have $\Phi(X)=2\tr(X)A$.

\emph{Step (ii).} Now we prove the general case. Since
$$\Phi(X)-J\Phi(X)^tJ^t=(A-C)X+X(B-D)-J[(A-C)X+X(B-D)]^tJ^t,$$
$$\Psi_{E\otimes F}(X)-J\Psi_{E\otimes F}(X)^tJ^t=\tr(EX)(F-JF^tJ^t),$$
the tuple $(A-C,B-D,C-A,D-B,E\otimes(F-JF^tJ^t))$ is $q$-admissible.
By Step (i), we have $\tr(A-C)=0$,
\begin{equation}\label{E:ABCD}
B-D=A-C,
\end{equation}
and
$$\Phi(X)-J\Phi(X)^tJ^t=2\tr(X)(A-C).$$ This also implies that the
traces of $A, B, C, D$ are the same. On the other hand, we have
\begin{align*}
\Phi(X)+J\Phi(X)^tJ^t&=(A+C)X+X(B+D)+J[(A+C)X+X(B+D)]^tJ^t\\
&=\tr((A+B+C+D)X)I\\
&=2\tr((A+D)X)I.
\end{align*}
Hence
\begin{equation}\label{E:2F}
\Phi(X)=\tr(X)(A-C)+\tr((A+D)X)I.
\end{equation}
There are two cases to consider.

\emph{Case (a).} $A-C$ and $I$ are linearly dependent. Then $A-C$ is
a scalar matrix. But we have $\tr(A)=\tr(C)$. Hence $C=A$. From
\eqref{E:ABCD}, we see that $D=B$. By \eqref{E:q}, we have
$\Phi=\Psi_{(A+B)\otimes I}$ and $\TT=\TT^q_{A,B}$.

\emph{Case (b).} $A-C$ and $I$ are linearly independent. Since
$\dim\im(\Phi)\leq1$, by \eqref{E:2F}, the dimension of the subspace
$$\{(\tr(X),\tr((A+D)X))\mid X\in\MM(2,\CC)\}$$ of $\CC^2$ is less
than or equal to $1$. This implies that $A+D$ is a scalar matrix. By
\eqref{E:ABCD}, $B+C$ is also a scalar matrix. Hence
$$D=A+D-A=\frac{1}{2}\tr(A+D)I-A=\tr(A)I-A=JA^tJ^t.$$
Similarly, we have $C=JB^tJ^t$. From \eqref{E:qdag}, we see that
$\Phi=\Psi_{I\otimes(A+B)}$ and $\TT=(\TT^q_{A,B})^\dag$.
\end{proof}

\begin{proof}[Proof of Proposition \ref{P:r} (2)]
By \eqref{E:FIJ}, the tuple $(A,-JBJ,-C,JDJ,-(JE)\otimes(FJ))$ is
$q$-admissible, which must be $\TT^q_{A,-JBJ}$ or
$(\TT^q_{A,-JBJ})^\dag$ by Proposition \ref{P:q} (1). This implies
that $\TT$ is equal to $\TT^r_{A,B}$ or $(\TT^r_{A,B})^\dag$.
\end{proof}

Finally we prove \ref{P:r} (3). We will make use of the
representations $\rho_1$ and $\rho_2$ of $\Lgl(3,\CC)$ in
$\MM_{\mathrm{skew}}(3,\CC)$ and $\CC^3$ defined in
\eqref{E:representation}.

\begin{proof}[Proof of Proposition \ref{P:r} (3)]
Denote $\Phi=\Phi_{A,B,C,D}^r$. Then $\dim\im(\Phi)\leq1$ and
$\Phi(X)_{ij}$ $(i,j\in\NN_3)$ are mutually linearly dependent. We
divide the proof into two steps.

\emph{Step (i).} First we assume that $C=A$, $D=B$. Then
$$\Phi(X)=AX+XB+(AX+XB)^t.$$ We prove that $A=B=0$.

Let $i,j\in\NN_3$, $i\neq j$. Let $k\in\NN_3$ with
$\{i,j,k\}=\NN_3$. Since
\begin{align*}
\Phi(X)_{ii}&=0X_{jk}+2A_{ij}X_{ji}+\cdots,\\
\Phi(X)_{ik}&=A_{ij}X_{jk}+\cdots,
\end{align*}
their linear dependence implies that $A_{ij}=0$. So $A$ is diagonal.
Similarly, $B$ is diagonal. Now we have six polynomials
\begin{align*}
\Phi(X)_{rr}=&2(A_{rr}+B_{rr})X_{rr}, && r\in\NN_3,\\
\Phi(X)_{rs}=&(A_{rr}+B_{ss})X_{rs}+(A_{ss}+B_{rr})X_{sr}, &&
(r,s)\in\{(1,2),(2,3),(3,1)\},
\end{align*}
whose mutual linear dependence forces that $A=-B\in\CC I$. But
$\tr(A)=\tr(B)$. So $A=B=0$.

\emph{Step (ii).} Now we prove the general case. Since
$$\Phi(X)+\Phi(X)^t=(A+C)X+X(B+D)+[(A+C)X+X(B+D)]^t,$$
$$\Psi_{E\otimes F}(X)+\Psi_{E\otimes F}(X)^t=\tr(EX)(F+F^t),$$ the tuple
$(A+C,B+D,A+C,B+D,E\otimes(F+F^t))$ is $r$-admissible. By Step (i),
we have $A+C=0$, $B+D=0$. So $$\Phi(X)=AX+XB-(AX+XB)^t.$$

Let $i,j\in\NN_3$, $i\neq j$. Let $k\in\NN_3$ with
$\{i,j,k\}=\NN_3$. Since
\begin{align*}
\Phi(X)_{ij}&=(A_{ij}-B_{ji})X_{jj}+\cdots,\\
\Phi(X)_{ik}&=0X_{jj}+A_{ij}X_{jk}-B_{ji}X_{kj}+\cdots,
\end{align*}
their linear dependence implies that
\begin{equation}\label{E:ss}
B_{ji}=A_{ij}.
\end{equation} We now prove that
\begin{equation}\label{E:Aii2}
A_{ii}-B_{ii}=A_{jj}-B_{jj}.
\end{equation}
If both $\Phi(X)_{ik}$ and $\Phi(X)_{jk}$ are identically zero, from
the expressions
\begin{align*}
\Phi(X)_{ik}&=(A_{ii}+B_{kk})X_{ik}-(A_{kk}+B_{ii})X_{ki}+\cdots,\\
\Phi(X)_{jk}&=(A_{jj}+B_{kk})X_{jk}-(A_{kk}+B_{jj})X_{kj}+\cdots,
\end{align*}
we get $A_{ii}+B_{kk}=A_{jj}+B_{kk}=A_{kk}+B_{ii}=A_{kk}+B_{jj}=0$,
which implies \eqref{E:Aii2}. If one of $\Phi(X)_{ik}$ and
$\Phi(X)_{jk}$, say $\Phi(X)_{ik}$, is not identically zero. Then we
have the linearly dependent polynomials
\begin{align*}
\Phi(X)_{ik}=&B_{jk}X_{ij}-A_{kj}X_{ji}+\cdots,\\
\Phi(X)_{ij}=&(A_{ii}+B_{jj})X_{ij}-(A_{jj}+B_{ii})X_{ji}+\cdots.
\end{align*}
Since $\Phi(X)_{ik}\not\equiv0$ and $B_{jk}=A_{kj}$, we must have
$A_{ii}+B_{jj}=A_{jj}+B_{ii}$, which also implies \eqref{E:Aii2}. By
\eqref{E:Aii2}, there exists $\alpha\in\CC$ such that
$B_{ii}=A_{ii}+\alpha$. But $\tr(A)=\tr(B)$. So we have
\begin{equation}\label{E:tt}
B_{ii}=A_{ii}.
\end{equation}
From \eqref{E:ss} and \eqref{E:tt}, we have $B=A^t$. Thus
$$\Phi(X)=A(X-X^t)+(X-X^t)A^t.$$ Denote
$\Phi_1=\Phi|_{\MM_{\mathrm{skew}}(3,\CC)}$. Then
$$\Phi_1(Y)=2(AY+YA^t)$$ for $Y\in\MM_{\mathrm{skew}}(3,\CC)$, and we
have $\dim\im(\Phi_1)\leq1$.

Now we consider the representations $\rho_1$ and $\rho_2$ of
$\Lgl(3,\CC)$ in $\MM_{\mathrm{skew}}(3,\CC)$ and $\CC^3$ defined in
\eqref{E:representation}. Note that $\Phi_1=2\rho_1(A)$. From the
proof of Lemma \ref{L:n=3}, we know that $\rho_1$ and $\rho_2$ are
equivalent. So
\begin{align*}
\rank(A-\tr(A)I)&=\rank(\tr(A)I-A^t)=\dim\im(\rho_2(A))\\
&=\dim\im(\rho_1(A))=\dim\im(\Phi_1)\leq1.
\end{align*}
Hence there exist $u,v\in\CC^3$ such that $A-\tr(A)I=uv^t$, i.e.,
$A=uv^t-\langle u,v\rangle I/2=\tau_{u,v}$. Now we have
$B=A^t=\tau_{v,u}$, $C=-A=-\tau_{u,v}$, and $D=-B=-\tau_{v,u}$. By
\eqref{E:tausigma}, we have $\Phi=\Psi_{\sigma_u\otimes\sigma_v}$.
Therefore $\TT=\TT_{u,v}$.
\end{proof}

\section{The main theorems}

Using the results about admissible tuples obtained in the previous
section, in this section we prove our main theorems (Theorems
\ref{T:nonhomogeneouspure}--\ref{T:nonhomogeneous} below). We first
prove a lemma, which is crucial for converting
Eq.~\eqref{E:eq-general} to matrix equations.

\begin{lemma}\label{L:matrixeq}
$\FF=(f_i)_{i=1}^6$ is a solution of Eq.~\eqref{E:eq-general} on $G$
if and only if
\begin{equation}\label{E:fourier2}
\tr[(\hat{f_1}(\pi)X+X\hat{f_2}(\pi))\pi(y)+(\hat{f_3}(\pi)X
+X\hat{f_4}(\pi))^t\bar{\pi}(y)]=\tr(\hat{f_5}(\pi)X)f_6(y)
\end{equation}
for all $y\in G$, $[\pi]\in\hat{G}$, and $X\in\MM(d_\pi,\CC)$.
\end{lemma}

\begin{proof}
Eq.~\eqref{E:eq-general} can be rewritten as
$$
R_yf_1+L_{y^{-1}}f_2+R_{y^{-1}}f_3+L_yf_4=f_6(y)f_5.
$$
Taking the Fourier transform, we see that this is equivalent to
$$
\pi(y)\hat{f_1}(\pi)+\hat{f_2}(\pi)\pi(y)+\pi(y)^{-1}\hat{f_3}(\pi)
+\hat{f_4}(\pi)\pi(y)^{-1}=f_6(y)\hat{f_5}(\pi)
$$
for all $[\pi]\in\hat{G}$. Then the lemma follows from the fact that
a matrix $A\in\MM(d_\pi,\CC)$ is equal to $0$ if and only if
$\tr(AX)=0$ for all $X\in\MM(d_\pi,\CC)$.
\end{proof}

In our first theorem we determine all pure normalized solutions of
Eq.~\eqref{E:eq-general}. We keep the notation from Examples
\ref{Ex:U1}--\ref{Ex:O3}.

\begin{theorem}\label{T:nonhomogeneouspure}
Let $[\pi]\in\hat{G}$, and let $\FF$ be a nontrivial pure normalized
solution of Eq.~\eqref{E:eq-general} on $G$ supported on $[[\pi]]$.
Denote $K=U(d_\pi)$, $O(d_\pi)$ or $Sp(d_\pi)$ according to the type
of $\pi$. Then $\FF=\FF^K\circ\pi$, where $\FF^K$ is a solution of
Eq.~\eqref{E:eq-general} on $K$, and the only possibilities of $K$
and $\FF^K$ are as follows:
\begin{itemize}
\item[(1)] $K=U(1)$ and
$\FF^{K}=\FF_{\varepsilon_1,\delta_1,\varepsilon_2,\delta_2}^{U(1)}$
for some $\varepsilon_1,\delta_1,\varepsilon_2,\delta_2\in\CC$;
\item[(2)] $K=O(2)$, $\FF^{K}=\FF^{O(2)}_{A,B}$ or
$(\FF^{O(2)}_{A,B})^\dag$ for some $A,B\in\MM(2,\CC)$ with
$\tr(A)=\tr(B)$;
\item[(3)] $K=SU(2)$, $\FF^{K}=\FF^{SU(2)}_{A,B}$ or
$(\FF^{SU(2)}_{A,B})^\dag$ for some $A,B\in\MM(2,\CC)$ with
$\tr(A)=\tr(B)$;
\item[(4)] $K=O(3)$ and $\FF^{K}=\FF_{u,v}^{O(3)}$ for some $u,v\in\CC^3$.
\end{itemize}
\end{theorem}

\begin{proof}
Since $\FF$ is normalized, we have
$$\tr(\hat{f_1}(\pi)-\hat{f_2}(\pi))=\tr(\hat{f_3}(\pi)-\hat{f_4}(\pi))=0.$$
According to the types of $\pi$ (c.f. Theorem \ref{T:types}), there
are three cases to consider.

\emph{Case (a).} $\pi$ is of complex type, i.e.,
$[\pi]\neq[\bar{\pi}]$. Applying Lemma \ref{L:matrixeq} to $\pi$ and
$\bar{\pi}$, we have
\begin{align*}
\hat{f_1}(\pi)X+X\hat{f_2}(\pi)&=\tr(\hat{f_5}(\pi)X)\hat{f_6}(\pi),\\
(\hat{f_3}(\pi)X+X\hat{f_4}(\pi))^t&=\tr(\hat{f_5}(\pi)X)\hat{f_6}(\bar{\pi}),\\
\hat{f_1}(\bar{\pi})X+X\hat{f_2}(\bar{\pi})&=\tr(\hat{f_5}(\bar{\pi})X)\hat{f_6}(\bar{\pi}),\\
(\hat{f_3}(\bar{\pi})X+X\hat{f_4}(\bar{\pi}))^t&=\tr(\hat{f_5}(\bar{\pi})X)\hat{f_6}(\pi)
\end{align*}
for all $X\in\MM(d_\pi,\CC)$. So the $3$-tuples
\begin{align*}
&(\hat{f_1}(\pi),\hat{f_2}(\pi),\hat{f_5}(\pi)\otimes\hat{f_6}(\pi)),
&&(\hat{f_3}(\pi),\hat{f_4}(\pi),\hat{f_5}(\pi)\otimes\hat{f_6}(\bar{\pi})^t),\\
&(\hat{f_1}(\bar{\pi}),\hat{f_2}(\bar{\pi}),\hat{f_5}(\bar{\pi})\otimes\hat{f_6}(\bar{\pi})),
&&(\hat{f_3}(\bar{\pi}),\hat{f_4}(\bar{\pi}),\hat{f_5}(\bar{\pi})\otimes\hat{f_6}(\pi)^t)
\end{align*}
are $c$-admissible. Since $\FF$ is nontrivial and supported on
$[[\pi]]$, these tuples can not be all trivial. By Proposition
\ref{P:c}, we have $d_\pi=1$, i.e., $K=U(1)$. Let
$\hat{f_5}(\pi)=\varepsilon_1$,
$\hat{f_5}(\bar{\pi})=\varepsilon_2$, $\hat{f_6}(\pi)=\delta_1$,
$\hat{f_6}(\bar{\pi})=\delta_2$. Then
\begin{align*}
\hat{f_1}(\pi)=\hat{f_2}(\pi)=\varepsilon_1\delta_1/2,\quad
\hat{f_1}(\bar{\pi})=\hat{f_2}(\bar{\pi})=\varepsilon_2\delta_2/2,\\
\hat{f_3}(\pi)=\hat{f_4}(\pi)=\varepsilon_1\delta_2/2,\quad
\hat{f_3}(\bar{\pi})=\hat{f_4}(\bar{\pi})=\varepsilon_2\delta_1/2.
\end{align*}
From Example \ref{Ex:U1} and the Fourier inversion formula, we see
that
$\FF=\FF_{\varepsilon_1,\delta_1,\varepsilon_2,\delta_2}^{U(1)}\circ\pi$.

\emph{Case (b).} $\pi$ is of real type, i.e., $\pi(G)\subseteq
O(d_\pi)$. Then $\bar{\pi}(y)=\pi(y)$ for all $y\in G$. By Lemma
\ref{L:matrixeq}, for all $X\in\MM(d_\pi,\CC)$ we have
$$\hat{f_1}(\pi)X+X\hat{f_2}(\pi)+(\hat{f_3}(\pi)X+X\hat{f_4}(\pi))^t
=\tr(\hat{f_5}(\pi)X)\hat{f_6}(\pi).$$ So the $5$-tuple
$$\TT_r=(\hat{f_1}(\pi),\hat{f_2}(\pi),\hat{f_3}(\pi),\hat{f_4}(\pi),\hat{f_5}(\pi)\otimes\hat{f_6}(\pi))$$
is $r$-admissible. Since $\FF$ is nontrivial and supported on
$[[\pi]]$, $\TT_r$ is nontrivial. By Proposition \ref{P:r}, we have
$d_\pi=1, 2$, or $3$, which correspond to the case of $K=O(1)$,
$O(2)$ or $O(3)$, respectively.

If $d_\pi=1$, then $\TT_r=\TT_{a,b}$
and $\FF=\FF_{a,b}^{O(1)}\circ\pi$ for some $a,b\in\CC$ (see Example
\ref{Ex:O1}). As mentioned in Example \ref{Ex:O1}, this case can be
absorbed into the case of $K=U(1)$.

If $d_\pi=2$, then
$\TT_r=\TT^r_{A,B}$ or $(\TT^r_{A,B})^\dag$, and
$\FF=\FF^{O(2)}_{A,B}\circ\pi$ or $(\FF^{O(2)}_{A,B})^\dag\circ\pi$
for some $A,B\in\MM(2,\CC)$ with $\tr(A)=\tr(B)$ (see Example \ref{Ex:O2}).

If $d_\pi=3$, then $\TT_r=\TT_{u,v}$ and
$\FF=\FF^{O(3)}_{u,v}\circ\pi$ for some $u,v\in\CC^3$ (see Example
\ref{Ex:O3}).

\emph{Case (c).}  $\pi$ is of quaternionic type, i.e., $d_\pi$ is
even and $\pi(G)\subseteq Sp(d_\pi)$. Then $\bar{\pi}(y)=J\pi(y)J^t$
for all $y\in G$. By Lemma \ref{L:matrixeq}, for all
$X\in\MM(d_\pi,\CC)$ we have
$$\hat{f_1}(\pi)X+X\hat{f_2}(\pi)+J(\hat{f_3}(\pi)X+X\hat{f_4}(\pi))^tJ^t=
\tr(\hat{f_5}(\pi)X)\hat{f_6}(\pi).$$ So the $5$-tuple
$$\TT_q=(\hat{f_1}(\pi),\hat{f_2}(\pi),\hat{f_3}(\pi),\hat{f_4}(\pi),\hat{f_5}(\pi)\otimes\hat{f_6}(\pi))$$
is $q$-admissible. As before, $\TT_q$ is nontrivial. By Proposition
\ref{P:q}, we have $d_\pi=2$ and $\TT_q=\TT^q_{A,B}$ or
$(\TT^q_{A,B})^\dag$. Hence $K=Sp(2)=SU(2)$, and
$\FF=\FF^{SU(2)}_{A,B}\circ\pi$ or
$(\FF^{SU(2)}_{A,B})^\dag\circ\pi$ for some $A,B\in\MM(2,\CC)$ with
$\tr(A)=\tr(B)$ (see Example \ref{Ex:SU2}).
\end{proof}

Our next theorem gives all pure normalized homogeneous solutions.

\begin{theorem}\label{T:homogeneouspure}
Under the same conditions as in Theorem \ref{T:nonhomogeneouspure},
if moreover $\FF$ is homogeneous, then the only possibilities of $K$
and $\FF^K$ are as follows:
\begin{itemize}
\item[(1)] $K=O(1)$ and $\FF^{K}=\FF_a^{O(1)}$ for some $a\in\CC$;
\item[(2)] $K=O(2)$ and $\FF^{K}=\FF^{O(2)}_A$
for some $A\in\MM(2,\CC)$ with $\tr(A)=0$;
\item[(3)] $K=SU(2)$ and $\FF^{K}=\FF^{SU(2)}_A$ for some $A\in\MM(2,\CC)$
with $\tr(A)=0$.
\end{itemize}
\end{theorem}

\begin{proof}
This follows directly from the proof of Theorem
\ref{T:nonhomogeneouspure} and the conditions for $\FF^K$ being
homogeneous given in Examples \ref{Ex:U1}--\ref{Ex:O3}.
\end{proof}

The next theorem characterizes the space of normalized homogeneous
solutions.

\begin{theorem}\label{T:homogeneous}
The Hilbert space of normalized homogeneous solutions of
Eq.~\eqref{E:eq-general} is spanned by pure normalized homogeneous
solutions.
\end{theorem}

\begin{proof}
Let $\FF=(f_i)_{i=1}^4$ be a normalized homogeneous solution of
Eq.~\eqref{E:eq-general} on $G$. It suffices to prove that $\FF$ is
the sum of some pure normalized homogeneous solutions. For
$\varpi\in[\hat{G}]$, let
$$f_i^{\varpi}(x)=\sum_{[\pi]\in\varpi}\tr(\hat{f_i}(\pi)\pi(x)), \quad 1\leq i\leq4.$$
Then
$$(f_i^{\varpi})\hat{}\,(\pi)=\begin{cases}\hat{f_i}(\pi), & [\pi]\in\varpi;\\0, & [\pi]\notin\varpi.\end{cases}$$
By Lemma \ref{L:matrixeq}, $\FF^\varpi=(f_i^\varpi)_{i=1}^4$ is a
pure normalized homogeneous solution of Eq.~\eqref{E:eq-general}
supported on $\varpi$, and we have
$\FF=\sum_{\varpi\in[\hat{G}]}\FF^\varpi$.
This proves the theorem.
\end{proof}

Finally we prove the theorem on the structure of the general
solution of Eq.~\eqref{E:eq-general} on $G$.

\begin{theorem}\label{T:nonhomogeneous}
Any solution of Eq.~\eqref{E:eq-general} on $G$
is of the form
$$\FF_0+\FF_h,$$
where $\FF_0$ is a
pure normalized solution and $\FF_h$ is a homogeneous solution.
\end{theorem}

\begin{proof}
Let $\FF=(f_i)_{i=1}^6$ be a solution of Eq.~\eqref{E:eq-general} on
$G$. Applying Lemma \ref{L:matrixeq} to $\FF$ and taking the Fourier
transform at the both sides of \eqref{E:fourier2}, we obtain that
$\supp((\tr(\hat{f_5}(\pi)X)f_6)\,\hat{}\,)\subseteq[[\pi]]$ for all
$[\pi]\in\hat{G}$ and $X\in\MM(d_{\pi},\CC)$. So if
$[\pi]\in\supp(\hat{f_5})$, then $\supp(\hat{f_6})\subseteq[[\pi]]$.
Hence there exists $\varpi_0\in[\hat{G}]$ such that
$\supp(\hat{f_5})\cup\supp(\hat{f_6})\subseteq\varpi_0$. Let
$$f_i^{\varpi_0}(x)=\sum_{[\pi]\in\varpi_0}\tr(\hat{f_i}(\pi)\pi(x)), \quad 1\leq i\leq4.$$
Then
$$(f_i^{\varpi_0})\hat{}\,(\pi)=\begin{cases}\hat{f_i}(\pi), & [\pi]\in\varpi_0;\\
0, & [\pi]\notin\varpi_0.\end{cases}$$ By Lemma \ref{L:matrixeq},
$\FF_0=(f_1^{\varpi_0},f_2^{\varpi_0},f_3^{\varpi_0},f_4^{\varpi_0},f_5\otimes
f_6)$ is a pure normalized solution supported on $\varpi_0$. So
$\FF_h=(f_i-f_i^{\varpi_0})_{i=1}^4$ is a homogeneous solution of
Eq.~\eqref{E:eq-general} on $G$, and we have $\FF=\FF_0+\FF_h$.
\end{proof}

Theorems \ref{T:nonhomogeneouspure}--\ref{T:nonhomogeneous}
provide a complete picture of the general solution
of Eq.~\eqref{E:eq-general} on the compact group $G$. They also
provide a method about how to construct all solutions. For a fixed
$G$, we first find all irreducible representations of $G$ into
$U(1)$, $O(2)$, $SU(2)$, and $O(3)$. Then using Theorem
\ref{T:nonhomogeneouspure}, we find all pure normalized solutions.
Theorem \ref{T:homogeneouspure} gives all pure normalized
homogeneous solutions. Here we should be careful that
representations into $O(1)$ provide nontrivial homogeneous
solutions. Theorem \ref{T:homogeneous} tells us that pure normalized
homogeneous solutions and solutions of the form $\FF_{c_1,c_2}$ span
the space of homogeneous solutions. Thus we determine all
homogeneous solutions. Finally by Theorem \ref{T:nonhomogeneous}, we
get the general solution by picking an arbitrary pure normalized
solution and take its sum with an arbitrary homogeneous solution. We
illustrate this by finding the general solution of
Eq.~\eqref{E:eq-general} on $SU(2)$.

\begin{example}[General Solution on $SU(2)$]
It is well known that for each positive integer $d$ there exists
exactly one $d$-dimensional irreducible representation of $SU(2)$
(see, e.g., \cite{BtD}). The $1$-dimensional one is
the trivial representation. So it is a representation into $O(1)$.
The $2$-dimensional one is the identity representation. The
$3$-dimensional one is the adjoint representation $\mathrm{Ad}$ in
the Lie algebra $\mathfrak{su}(2)$ of $SU(2)$, which can be viewed
as a representation into $O(3)$. As the $1$-dimensional
representation is into $O(1)$, when applying Theorem
\ref{T:nonhomogeneouspure} (1), we can use Example \ref{Ex:O1}.
Indeed, as the $1$-dimensional representation is trivial, the pure
normalized solutions obtained from Theorem
\ref{T:nonhomogeneouspure} (1) are constant solutions. They are of
the form
\begin{equation}\label{E:5.2}
f_1\equiv f_2\equiv a/2, \quad f_3\equiv f_4\equiv b/2, \quad
f_5\otimes f_6\equiv a+b
\end{equation}
for some $a,b\in\CC$. The pure normalized solutions obtained by
applying Theorem \ref{T:nonhomogeneouspure} (3)--(4) to the identity
representation and the adjoint representation are
$\FF^{SU(2)}_{A,B}$, $(\FF^{SU(2)}_{A,B})^\dag$, and
$\FF^{O(3)}_{u,v}\circ\mathrm{Ad}$. Thus we get all pure normalized
solutions of Eq.~\eqref{E:eq-general} on $SU(2)$. Now applying
Theorem \ref{T:homogeneouspure}, we obtain all pure normalized
homogeneous solutions. They are $f_1\equiv f_2\equiv -f_3\equiv
-f_4\equiv\mathrm{const}$ and $\FF^{SU(2)}_{A'}$. By Theorem
\ref{T:homogeneous}, all homogeneous solutions of
Eq.~\eqref{E:eq-general} on $SU(2)$ are of the form
\begin{equation}\label{E:5.6}
\begin{cases}
f_1(x)=\tr(A'x)+c_1(x)+\alpha,\\
f_2(x)=-\tr(A'x)-c_1(x)+\alpha,\\
f_3(x)=\tr(A'x)+c_2(x)-\alpha,\\
f_4(x)=-\tr(A'x)-c_2(x)-\alpha,
\end{cases}
\end{equation}
where $A'\in\MM(2,\CC)$, $c_1,c_2\in L^2_c(G)$,
$\alpha\in\CC$. Finally, by Theorem \ref{T:nonhomogeneous}, the
general solution of Eq.~\eqref{E:eq-general} on $SU(2)$ is given by
$\FF_0+\FF_h$, where $\FF_0\in\{\eqref{E:5.2}, \FF^{SU(2)}_{A,B},
(\FF^{SU(2)}_{A,B})^\dag, \FF^{O(3)}_{u,v}\circ\mathrm{Ad}\}$ and
$\FF_h$ is given by \eqref{E:5.6}.
\end{example}

\section{Applications}

In this section, we consider some functional equations on compact
groups which are special cases of Eq.~\eqref{E:eq-general}. In
particular, we solve the Wilson equation and the d'Alembert long
equation on compact groups. We also recover the general solution of
the d'Alembert equation that was obtained in \cite{Davison1,
Yang06}.

We first consider the equation
\begin{eqnarray}\label{E:Wilson_gen1}
f(xy)+g(xy^{-1})=h(x)k(y),
\end{eqnarray}
where $f,g,h,k:G\to \CC$ are the unknowns. It is clear that
Eq.~\eqref{E:Wilson_gen1} corresponds to the special case of
Eq.~\eqref{E:eq-general} where $f_2\equiv f_4\equiv0$. We denote a
solution of Eq.~\eqref{E:Wilson_gen1} by $\FF=(f,g,h\otimes k)$, and
say that it is \emph{homogeneous} if $h\otimes k\equiv0$. If
$\FF=(f,g,h\otimes k)$ is a solution and $\FF'=(f',g',0)$ is a
homogeneous solution of Eq.~\eqref{E:Wilson_gen1}, then
$\FF+\FF'=(f+f',g+g',h\otimes k)$ is also a solution of
Eq.~\eqref{E:Wilson_gen1}. We first construct some homogeneous
solutions of Eq.~\eqref{E:Wilson_gen1}.

\begin{example}\label{Ex:6.4}
Let $\pi:G\rightarrow O(1)$ be a homomorphism, and let $a\in\CC$. We
view $\pi$ as a function on $G$. Then
$$\FF_{\pi,a}=(a\pi,-a\pi,0)$$
is a homogeneous solution
of Eq.~\eqref{E:Wilson_gen1} on $G$. More generally, if
$\pi_j:G\rightarrow O(1)$ are distinct homomorphisms and $a_j\in\CC$
($j=1,2,\ldots$), then
$$\sum_{j\geq1}\FF_{\pi_j,a_j}=(\sum_{j\geq1}a_j\pi_j,-\sum_{j\geq1}a_j\pi_j,0)$$
is a homogeneous solution, provided that
$\sum_{j\geq1}|a_j|^2<\infty$.
\end{example}

Now we construct some solutions of Eq.~\eqref{E:Wilson_gen1} on
$U(1)$, $O(2)$, and $SU(2)$.

\begin{example}\label{Ex:6.1}
Let $G=U(1)$. For
$\varepsilon_1,\delta_1,\varepsilon_2,\delta_2\in\CC$, define
$$\begin{cases}
f(x)=\varepsilon_1\delta_1x+\varepsilon_2\delta_2\bar{x},\\
g(x)=\varepsilon_1\delta_2x+\varepsilon_2\delta_1\bar{x}, \\
h\otimes
k(x,y)=(\varepsilon_1x+\varepsilon_2\bar{x})(\delta_1y+\delta_2\bar{y}),
\end{cases} \quad x,y\in U(1).$$
It is easy to check that $(f,g,h\otimes k)$ is a solution of
Eq.~\eqref{E:Wilson_gen1} on $U(1)$.
\end{example}

\begin{example}\label{Ex:6.2}
Let $G=O(2)$. For $P\in\MM(2,\CC)$, define
$$\begin{cases}
f(x)=-g(x)=\tr(Px), \\
h\otimes k(x,y)=-\tr(JPx)\tr(Jy),
\end{cases} \quad x,y\in O(2).$$
Then $(f,g,h\otimes k)$ is a solution of Eq.~\eqref{E:Wilson_gen1}
on $O(2)$.
\end{example}

\begin{example}\label{Ex:6.3}
Let $G=SU(2)$. For $P\in\MM(2,\CC)$, define
$$
\begin{cases}
f(x)=g(x)=\tr(Px), \\
h\otimes k(x,y)=\tr(Px)\tr(y),
\end{cases} \quad x,y\in SU(2).$$
Then $(f,g,h\otimes k)$ is a solution  of Eq.~\eqref{E:Wilson_gen1}
on $SU(2)$.
\end{example}

We leave the verification of the above examples to the reader. The
following result claims that the above examples are the building
blocks of the general solution of Eq.~\eqref{E:Wilson_gen1} on $G$.

\begin{theorem}\label{T:gen1}
Any solution of Eq.~\eqref{E:Wilson_gen1} on $G$ is of the form
$$\FF\circ\pi+\sum_{j\geq1}\FF_{\pi_j,a_j},$$
where $\pi:G\rightarrow K$ is an irreducible representation with
$K=U(1)$, $O(2)$ or $SU(2)$, $\FF$ is a solution of
Eq.~\eqref{E:Wilson_gen1} on $K$ as in Examples
\ref{Ex:6.1}--\ref{Ex:6.3}, and $\sum_{j\geq1}\FF_{\pi_j,a_j}$ as in
Example \ref{Ex:6.4}.
\end{theorem}

\begin{proof}
Let $(f,g,h\otimes k)$ be a solution of Eq.~\eqref{E:Wilson_gen1}.
Then $(f,0,g,0,h\otimes k)$ is a solution of
Eq.~\eqref{E:eq-general}. By Theorems
\ref{T:nonhomogeneouspure}--\ref{T:nonhomogeneous}, there exist
$c_1,c_2\in L^2_c(G)$ and irreducible representations
$\pi_j:G\rightarrow K_j$ ($j\geq0$) with $[[\pi_j]]$'s distinct,
such that
\begin{equation}\label{E:rhs}
(f,0,g,0,h\otimes
k)=\FF_{c_1,c_2}+\sum_{j\geq0}\FF^{K_j}\circ\pi_j,
\end{equation}
where
$\FF^{K_0}=(f_1^{K_0},f_2^{K_0},f_3^{K_0},f_4^{K_0},f_5^{K_0}\otimes
f_6^{K_0})$ is a solution of Eq.~\eqref{E:eq-general} on $K_0$,
$\FF^{K_j}=(f_1^{K_j},f_2^{K_j},f_3^{K_j},f_4^{K_j})$ ($j\geq1$) is
a homogeneous solution of Eq.~\eqref{E:eq-general} on $K_j$, and the
only possibilities of $K_j$, $\pi_j$, and $\FF^{K_j}$ are given in
Theorems \ref{T:nonhomogeneouspure} and \ref{T:homogeneouspure}.
Note that this implies
$$c_1=\sum_{j\geq0}f_2^{K_j}\circ\pi_j, \quad
c_2=\sum_{j\geq0}f_4^{K_j}\circ\pi_j$$ and
\begin{equation}\label{E:6.2}
f=\sum_{j\geq0}(f_1^{K_j}+f_2^{K_j})\circ\pi_j, \quad
g=\sum_{j\geq0}(f_3^{K_j}+f_4^{K_j})\circ\pi_j.
\end{equation}
Without loss of generality, we may assume that each $\FF^{K_j}$ is a
nontrivial solution.

We first prove that $K_0\neq O(3)$. Suppose $K_0=O(3)$. Then
$\FF^{K_0}=\FF_{u,v}^{O(3)}$ for some $u,v\in\CC^3$. Since
$\FF^{K_j}\circ\pi_j$ is a pure solution of Eq.~\eqref{E:eq-general}
on $G$ supported on $[[\pi_j]]$ for any $j\geq0$ and $[[\pi_j]]$ are
distinct, we have $(f_2^{K_j}\circ\pi_j)\hat{}\,(\pi_0)=0$ if
$j\geq1$. Hence
$$\hat{c}_1(\pi_0)=\sum_{j\geq0}(f_2^{K_j}\circ\pi_j)\hat{}\,(\pi_0)
=(f_2^{K_0}\circ\pi_0)\hat{}\,(\pi_0)=\tau_{v,u},$$ where
$\tau_{v,u}$ is as in Lemma \ref{L:n=3}. Since $c_1$ is a central
function, $vu^t=\hat{c}_1(\pi_0)+\langle u,v\rangle I_3/2$ is a
scalar matrix. This implies that $vu^t=0$, i.e., $u=0$ or $v=0$.
Hence $\FF^{K_0}$ is the trivial solution, a contradiction.

Now we prove that if $K_j=O(2)$, then $j=0$ and
$\FF^{K_0}=\FF^{O(2)}_{A,\frac{1}{2}\tr(A)I}$ for some
$A\in\MM(2,\CC)$. We know that if $K_j=O(2)$, then
$\FF^{K_j}=\FF^{O(2)}_{A,B}$ or $(\FF^{O(2)}_{A,B})^\dag$ for some
$A,B\in\MM(2,\CC)$ with $\tr(A)=\tr(B)$, and $B=A^t$ with $\tr(A)=0$
if $j\geq1$. If $\FF^{K_j}=\FF^{O(2)}_{A,B}$, similar to the above
proof, we obtain that $B=\hat{c}_1(\pi_j)$ is a scalar matrix. So
$B=\tr(A)I/2$. If $j\geq1$, then $A=B=0$, conflicting with the
assumption that $\FF^{K_j}$ is nontrivial. Hence $j=0$. If
$\FF^{K_j}=(\FF^{O(2)}_{A,B})^\dag$, then similarly
$B=\hat{c}_1(\pi_j)$ and $-A^t=\hat{c}_2(\pi_j)$ are scalar
matrices. So $A=B=\lambda I$ for some $\lambda\in\CC$. By Remark
\ref{R:3.1}, this case can be absorbed into the former case. Note
that if we set $P=A+\tr(A)I/2$, then we have
\begin{equation}\label{E:6.3}
\begin{cases}f_1^{K_0}(x)+f_2^{K_0}(x)=-(f_3^{K_0}(x)+f_4^{K_0}(x))=\tr(Px),\\
f_5^{K_0}\otimes f_6^{K_0}(x,y)=-\tr(JPx)\tr(Jy),\end{cases} x,y\in
O(2).
\end{equation}

A similar argument shows that if $K_j=SU(2)$, then $j=0$ and
$\FF^{K_0}=\FF^{SU(2)}_{A,\frac{1}{2}\tr(A)I}$ for some
$A\in\MM(2,\CC)$. In this case if we set $P=A+\tr(A)I/2$, then we
have
\begin{equation}\label{E:6.4}
\begin{cases}f_1^{K_0}(x)+f_2^{K_0}(x)=f_3^{K_0}(x)+f_4^{K_0}(x)=\tr(Px)\\
f_5^{K_0}\otimes f_6^{K_0}(x,y)=\tr(Px)\tr(y),\end{cases} x,y\in
SU(2).
\end{equation}

The above proofs also imply that if $j\geq1$, then $K_j=O(1)$ and
$\FF^{K_j}=\FF^{O(1)}_{a_j}$ for some $a_j\in\CC$. In this case we
have
\begin{equation}\label{E:6.5}
f_1^{K_j}(x)+f_2^{K_j}(x)=-(f_3^{K_j}(x)+f_4^{K_j}(x))=a_jx, \quad
x\in O(1).
\end{equation}

Now we know that there are three possibilities for $K_0$, i.e.,
$K_0=U(1)$, $O(2)$, or $SU(2)$. In each case, it is easy to see from
\eqref{E:6.2}--\eqref{E:6.5} that
$$(f,g,h\otimes k)=\FF\circ\pi_0+\sum_{j}\FF_{\pi_j,a_j},$$
where $\FF$ is
a solution of Eq.~\eqref{E:Wilson_gen1} on $K_0$ as
in Examples \ref{Ex:6.1}--\ref{Ex:6.3}. The proof of the theorem is
completed by setting $K=K_0$ and $\pi=\pi_0$.
\end{proof}

Now we consider the special case of Eq.~\eqref{E:Wilson_gen1} where
$f\equiv g$.

\begin{theorem}\label{T:gen2}
The general solution of the equation
\begin{eqnarray}\label{E:Wilson_gen2}
f(xy)+f(xy^{-1})=h(x)k(y)
\end{eqnarray}
is
$$\begin{cases}
f(x)=\tr(P\pi(x)), \\
h\otimes k(x,y)=\tr(P\pi(x))\tr(\pi(y)),\end{cases}$$ where
$\pi:G\rightarrow SU(2)$ is a homomorphism and $P\in\MM(2,\CC)$.
\end{theorem}

\begin{proof}
Clearly, the general solution of Eq.~\eqref{E:Wilson_gen2}
corresponds to the solutions of Eq.~\eqref{E:Wilson_gen1} for which
$f\equiv g$. By Theorem \ref{T:gen1}, the functions $f$ and $g$ in a
solution of Eq.~\eqref{E:Wilson_gen1} has the form
$$f=f^K\circ\pi+\sum_{j\geq1}a_j\pi_j, \quad g=g^K\circ\pi-\sum_{j\geq1}a_j\pi_j,$$
where $K=U(1)$, $O(2)$ or $SU(2)$, $\pi:G\rightarrow K$ and
$\pi_j:G\rightarrow O(1)$ are distinct irreducible representations,
$f^K$ and $g^K$ are functions on $K$ as in Examples
\ref{Ex:6.1}--\ref{Ex:6.3}. Applying the Fourier transform, it is
easy to see that $f\equiv g$ if and only if $f^K\equiv g^K$ and
$a_j=0$. Restricting our attention to nontrivial solutions, we can
see that either $K=U(1)$ and $\delta_1=\delta_2$ (in the notation of
Example \ref{Ex:6.1}), or $K=SU(2)$. If $K=SU(2)$ we have reached
the conclusion of the theorem. If $K=U(1)$ and
$\delta_1=\delta_2=:\delta$, then the homomorphism
$x\mapsto\diag(\pi(x),\bar{\pi}(x))\in SU(2)$ and
$P=\diag(\varepsilon_1\delta,\varepsilon_2\delta)$ satisfy our
requirements.
\end{proof}

The following corollaries are straightforward from Theorem
\ref{T:gen2}.

\begin{corollary}\label{C:Wilson}
Any nontrivial solution of the Wilson equation \eqref{E:Wilson} is
of the form
$$\begin{cases}f(x)=\tr(P\pi(x)), \\ g(x)=\frac{1}{2}\tr(\pi(x)),\end{cases}$$ where
$\pi:G\rightarrow SU(2)$ is a homomorphism and $P\in\MM(2,\CC)$.
\end{corollary}

\begin{corollary}\label{C:vWilson}
Any nontrivial solution of the equation
\begin{equation}\label{E:vWilson}
f(xy)+f(xy^{-1})=2g(x)f(y)
\end{equation}
is of the form
$$\begin{cases}f(x)=a\tr(\pi(x)), \\ g(x)=\frac{1}{2}\tr(\pi(x)),\end{cases}$$ where
$\pi:G\rightarrow SU(2)$ is a homomorphism and $a\in\CC$.
\end{corollary}

\begin{corollary}\label{C:cosine}
Any nontrivial solution of the d'Alembert equation
\eqref{E:D'Alembert} is of the form
$$f(x)=\frac{1}{2}\tr(\pi(x)),$$ where $\pi:G\rightarrow SU(2)$ is
a homomorphism.
\end{corollary}

Indeed, to prove Corollaries \ref{C:Wilson}--\ref{C:cosine}, it
suffices to examine the solutions of Eq.~\eqref{E:Wilson_gen2}
satisfying $h\equiv 2f$, $k\equiv 2f$, and $h\equiv 2k\equiv 2f$,
respectively.

\medskip

Now we apply the results in the previous section to another type of
equations.

\begin{theorem}\label{T:long1}
Let $(f,h\otimes k)$ be a solution of the equation
\begin{equation}\label{E:longWilson}
f(xy)+f(xy^{-1})+f(yx)+f(y^{-1}x)=h(x)k(y).
\end{equation}
Then either there exist an irreducible
representation $\pi:G\rightarrow O(2)$ and $a\in\CC$ such that
$$\begin{cases}f(x)=a\tr(J\pi(x)), \\ h\otimes k(x,y)=2a\tr(J\pi(x))\tr(\pi(y)),\end{cases}$$
or there exist a representation $\pi:G\rightarrow SU(2)$ and
$A\in\MM(2,\CC)$ such that
$$\begin{cases}f(x)=\tr(A\pi(x)), \\ h\otimes
k(x,y)=2\tr(A\pi(x))\tr(\pi(y)).\end{cases}$$
\end{theorem}

\begin{proof}
It suffices to consider the solutions of Eq.~\eqref{E:eq-general}
satisfying $f_1\equiv f_2\equiv f_3\equiv f_4$. We write the general
solution of Eq.~\eqref{E:eq-general} as the right hand side of
\eqref{E:rhs}. In particular, we have
$$f_1=c_1+\sum_{j\geq0}f_1^{K_j}\circ\pi_j, \quad
f_2=-c_1+\sum_{j\geq0}f_2^{K_j}\circ\pi_j.$$ So $f_1\equiv f_2$
implies that
$$2c_1+\sum_{j\geq0}(f_1^{K_j}-f_2^{K_j})\circ\pi_j=0.$$ But the two
summands above belong to $L_c^2(G)$ and $L_c^2(G)^\bot$,
respectively. So we have $c_1\equiv0$. Similarly, we have
$c_2\equiv0$. By considering the Fourier transform, it is easy to
see that $f_1^{K_j}\equiv f_2^{K_j}\equiv f_3^{K_j}\equiv f_4^{K_j}$
for any $j\geq0$. Now one can verify that $\FF^{K_j}$ is
trivial if $j\geq1$, and $K_0$, $\pi_0$, and
$\FF^{K_0}$ take one of the following forms:
\begin{itemize}
\item[(1)] $K_0=U(1)$ and
$\FF^{K_0}=\FF_{\varepsilon_1,\delta,\varepsilon_2,\delta}^{U(1)}$
for some $\varepsilon_1,\varepsilon_2,\delta\in\CC$;
\item[(2)] $K_0=O(2)$ and $\FF^{K_0}=(\FF^{O(2)}_{aJ,aJ})^\dag$
for some $a\in\CC$;
\item[(3)] $K_0=SU(2)$ and $\FF^{K_0}=\FF^{SU(2)}_{A,A}$ for some
$A\in\MM(2,\CC)$.
\end{itemize}
The last two cases obviously satisfy the conclusion of the theorem.
For the first case, it suffices to set
$\pi(x)=\diag(\pi_0(x),\bar{\pi}_0(x))\in SU(2)$ and
$A=\diag(\varepsilon_1\delta,\varepsilon_2\delta)$.
\end{proof}

Similar to Corollaries \ref{C:Wilson}--\ref{C:cosine}, we have the
following corollaries.

\begin{corollary}\label{C:vlongWilson1}
Let $(f,g)$ be a nontrivial solution of the equation
\begin{equation}\label{E:vlongWilson1}
f(xy)+f(xy^{-1})+f(yx)+f(y^{-1}x)=4f(x)g(y).
\end{equation}
Then either there exist an irreducible
representation $\pi:G\rightarrow O(2)$ and $a\in\CC$ such that
$$\begin{cases}f(x)=a\tr(J\pi(x)), \\ g(x)=\frac{1}{2}\tr(\pi(x)),\end{cases}$$
or there exist a representation $\pi:G\rightarrow SU(2)$ and
$A\in\MM(2,\CC)$ such that
$$\begin{cases}f(x)=\tr(A\pi(x)), \\
g(x)=\frac{1}{2}\tr(\pi(x)).\end{cases}$$
\end{corollary}

\begin{corollary}\label{C:vlongWilson2}
Any nontrivial solution of the equation
\begin{equation}\label{E:vlongWilson2}
f(xy)+f(xy^{-1})+f(yx)+f(y^{-1}x)=4g(x)f(y)
\end{equation}
is of the form
$$\begin{cases}f(x)=a\tr(\pi(x)), \\ g(x)=\frac{1}{2}\tr(\pi(x)),\end{cases}$$
where $\pi:G\to SU(2)$ is a homomorphism and $a\in\CC$.
\end{corollary}

\begin{corollary}\label{C:longcosine}
Any nontrivial solution of the d'Alembert long equation
\eqref{E:long} is of the form
$$f(x)=\frac{1}{2}\tr(\pi(x)),$$ where $\pi: G\to SU(2)$ is a
homomorphism.
\end{corollary}

From Corollary \ref{C:longcosine}, we see that the solutions of the
d'Alembert long equation \eqref{E:long} and the d'Alembert equation
\eqref{E:D'Alembert} are the same. The similar result for step $2$
nilpotent groups was proved in \cite{Stet2}.

\medskip

The factorization property of the d'Alembert equation on compact
groups was studied in \cite{Davison1,Davison2,Yang06,Yangthesis}. To
conclude this section, we summarize the same property of the above
equations as follows.

\begin{corollary}
The following factorization properties hold.
\begin{itemize}
\item[(1)]
All nontrivial solutions of
Eqs.~\eqref{E:Wilson_gen2} and \eqref{E:vlongWilson2} on a compact group factor through $SU(2)$.
\item[(2)] All nontrivial solutions of Eq.~\eqref{E:longWilson} on a compact group
factor through $O(2)$ or $SU(2)$.
\end{itemize}
\end{corollary}

As a simple consequence, all nontrivial solutions of every special
case of Eqs.~\eqref{E:Wilson_gen2} and \eqref{E:vlongWilson2}, in
particular, the Wilson equation and the d'Alembert long equation,
factor through $SU(2)$.

\end{document}